# Multi-beam Beamforming in RIS-aided MIMO Subject to Reradiation Mask Constraints – Optimization and Machine Learning Design

Shumin Wang, Hajar El Hassani, Marco Di Renzo, *Fellow, IEEE*, and Marios Poulakis, *Senior Member, IEEE*

*Abstract*—Reconfigurable intelligent surfaces (RISs) are an emerging technology for improving spectral efficiency and reducing power consumption in future wireless systems. This paper investigates the joint design of the transmit precoding matrices and the RIS phase shift vector in a multi-user RIS-aided multiple-input multiple-output (MIMO) communication system. We formulate a max-min optimization problem to maximize the minimum achievable rate while considering transmit power and reradiation mask constraints. The achievable rate is simplified using the Arimoto-Blahut algorithm, and the problem is broken into quadratic programs with quadratic constraints (QPQC) sub-problems using an alternating optimization approach. To improve efficiency, we develop a model-based neural network optimization that utilizes the one-hot encoding for the angles of incidence and reflection. We address practical RIS limitations by using a greedy search algorithm to solve the optimization problem for discrete phase shifts. Simulation results demonstrate that the proposed methods effectively shape the multi-beam radiation pattern towards desired directions while satisfying reradiation mask constraints. The neural network design reduces the execution time, and the discrete phase shift scheme performs well with a small reduction of the beamforming gain by using only four phase shift levels.

## I. Introduction

Manuscript received Dec. 31, 2024; revised May 5, 2025, and July 11, 2025.

Part of this work was presented at IEEE ISWCS 2024 [1].

This research work was supported in part by the European Union through the Horizon 2020 project MetaWireless under grant agreement number 956256. The work of M. Di Renzo was supported in part by the European Union through the Horizon Europe project COVER under grant agreement number 101086228, the Horizon Europe project UNITE under grant agreement number 101129618, the Horizon Europe project INSTINCT under grant agreement number 101139161, and the Horizon Europe project TWIN6G under grant agreement number 101182794, as well as by the Agence Nationale de la Recherche (ANR) through the France 2030 project ANR-PEPR Networks of the Future under grant agreement NF-YACARI 22-PEFT-0005, and by the CHIST-ERA project PASSIONATE under grant agreements CHIST-ERA-22-WAI-04 and ANR-23-CHR4-0003-01.

S. Wang and M. Di Renzo are with Université Paris-Saclay, CNRS, CentraleSupélec, Laboratoire des Signaux et Systèmes, 3 Rue Joliot-Curie, 91192 Gif-sur-Yvette, France. (marco.direnzo@centralesupelec.fr). M Di renzo is also with King's College London, Centre for Telecommunications Research – Department of Engineering, WC2R 2LS London, United Kingdom (marco.di_renzo@kcl.ac.uk).

H. El Hassani is with ETIS, UMR 8051, CY Cergy Paris Université, ENSEA, CNRS, F-95000 Cergy, France. (hajar.el-hassani@ensea.fr).

M. Poulakis is with Huawei Technologies Sweden AB, Skalholtsgatan 9-11, 164 40 Kista, Sweden (marios.poulakis@huawei.com).

IN recent years, reconfigurable intelligent surfaces (RISs) have become a key technology for the advancement of sixth-generation (6G) networks due to their ability to significantly improve wireless communication performance at low hardware costs and power consumption [2]–[5]. Generally speaking, an RIS is a planar engineered surface composed of a large number of reconfigurable scattering elements, each capable of realizing advanced wave transformations by manipulating the amplitudes, phases and polarizations of incident waves [6]. This allows for optimal control over reflection, refraction, and other wave behaviors based on users' locations.

These appealing properties of RIS have motivated its application in various wireless scenarios, such as extending coverage to users located in dead zones and cell edges [7], [8], enabling massive device-to-device (D2D) communications [9], and supporting unmanned aerial vehicle (UAV) communications [10]. An RIS also enhances physical layer security [11] and facilitates simultaneous wireless information and power transfer (SWIPT) [12] and energy-efficient transmission schemes [13]. Recently, RIS has also been applied to integrated sensing and communication (ISAC) [14] and optical wireless systems [15], with system-level simulations demonstrating substantial benefits in large-scale network deployments [16].

### A. Related Works

Given the potential benefits of RIS, substantial research efforts have focused on developing advanced beamforming designs to enhance the performance of RIS-aided communication systems, as highlighted in [17]. For instance, the authors of [5] proposed an approach to enhance the spectrum and energy efficiency of multiple-input single-output (MISO) systems by jointly optimizing the transmit beamforming and the nearly passive beamforming at the RIS [1]. To further refine the beamforming design, the authors of [18] addressed practical hardware impairments (HWIs) and the dynamic noise introduced by active RIS.

In multiple-input multiple-output (MIMO) systems, numerous studies have explored RIS beamforming design [19]–[24]. For example, [19] and [20] maximized the capacity 2of an RIS-aided single-user (SU)-MIMO communication system by jointly optimizing the RIS reflection coefficients and the transmit covariance matrix using alternating optimization (AO). The authors of [21] introduced a low-complexity, RIS-partitioning-based scalable beamforming design to enhance the performance of large-scale MIMO systems. To support downlink transmission to cell-edge users, [22] proposed an algorithm that jointly optimizes the active precoder matrices at the base stations (BSs) and the phase shifts at the RIS, improving the strength of the desired signal and reducing cochannel interference. In [23], an AO-based method was introduced to design the source precoders and the RIS phase shift matrix for an RIS-aided full-duplex MIMO two-way communication system, aiming to maximize the system sum-rate. Practical RIS constraints, such as discrete phase shift levels, were also discussed in [24].

---

[1]Nearly passive beamforming refers to RIS operation where the phase of incoming signals is adjusted using low-power components such as PIN diodes or varactors. Unlike active systems that generate and amplify signals, RIS consumes minimal power for tuning, which makes it distinct from both active and fully passive systems.

Recent works have explored the use of RIS to improve multi-user (MU) transmission and beam coverage under various system constraints. By carefully adjusting the phase shifts of its reflective elements, RIS can align the reflected signals to create multiple beams, where each beam is directed towards a specific direction or user. This makes RIS a promising solution for supporting MU transmission. For example, [25] addressed the sum-rate maximization in downlink MU-MIMO systems with both continuous and discrete RIS phase shifts, using the weighted minimum mean square error (WMMSE) approach and Schnorr Euchner sphere decoding (SESD) algorithm. Using iterative optimization, [26] and [27] proposed hybrid beamforming schemes that combine digital beamforming at the BS with discrete analog beamforming at the RIS, showing that high sum-rate performance can be achieved even with limited phase resolution and partial channel knowledge. Similarly, [28] exploited the MIMO broadcast channel, multiple-access channel (BC-MAC) duality to maximize the sum-rate of an RIS-assisted MU-MIMO system by alternately optimizing the covariance matrices of the users and the RIS phase shifts. The authors of [29] introduced a method based on constructive interference, to jointly design RIS phase shifts, which aims to reduce symbol error rates and improve signal alignment at each user. In a different context, [30] studied RIS-assisted systems with both cooperative and non-cooperative BSs, where RIS configuration, active beamforming, power allocation and user association at BSs are jointly optimized to maximize the sum-rate for users served by cooperative BSs, while limiting RIS interference to non-cooperative users. For beam control, [31] presented an analytical multi-beam reconfiguration method that enables independent control of multiple beams without relying on heavy optimization. The method shows significant gains in system throughput in both indoor and outdoor scenarios. These works highlight the potential of RIS to support efficient MU transmission and beam coverage in practical systems.

Beamforming design has also been explored in other system setups. In [32], a semidefinite relaxation (SDR)-based AO algorithm was proposed for the joint design of active beamforming and phase shifts in an RIS-aided radar-communication (Radcom) system, demonstrating the benefits of RIS in enhancing Radcom performance. In [33], an energy-efficient SWIPT-enabled RIS-assisted MIMO system was investigated, optimizing both the active beamforming matrix at the BS and the nearly passive beamforming matrix at the RIS. Several studies have explored RIS capabilities in generating and steering multiple beams towards various user locations. Specifically, [34] and [35] examined multiple beam designs through beam training techniques. [34] developed a deep convolutional neural network to compute the RIS phase shift matrix, enabling multi-beam steering according to the desired beam pattern. Additionally, [35] demonstrated the application of deep neural networks for beam training in RIS-aided millimeter-wave (mmWave) massive MIMO systems. Despite the promising performance of these neural network-based methods in multi-beam design, they have drawbacks, such as requiring training datasets and corresponding output labels, which increase pilot overhead and computational complexity.

As an alternative, model-driven neural networks are becoming increasingly used to solve challenging optimization problems in wireless communications [36], including deep unfolding [37]–[39] and graph neural networks [40]. These methods have been applied to RIS-aided communications to effectively address non-convex optimization problems. For instance, the authors of [41] developed a deep unfolding model for channel estimation, achieving better performance with lower training overhead and computational complexity compared to the least square (LS) method. The authors of [42] introduced a deep denoising neural network combined with compressive sensing techniques for channel estimation in millimeter wave channels. In [43], a model-based deep unfolding technique was proposed to balance the tradeoff between communication rate and sensing accuracy with reduced complexity. Furthermore, [44] proposed a graph neural network-based approach to leverage the graph topology inherent in the optimization problem.

In the communication literature, optimization problems typically focus on shaping the RIS beampattern towards specified directions of interest, ignoring potential reradiated beams in other directions [3]. However, recent research works have highlighted the importance of applying reradiation constraints across the entire reradiation pattern of the RIS [3], [45]–[49]. Without proper control, an RIS can radiate power towards unintended directions where it may cause interference to other users or systems. Hence, it is necessary to impose specified reradiation mask constraints to better manage the interference and make the system more practical and compliant with real-world requirements and regulations. Motivated by these considerations, we recently investigated the design of a two-user RIS-aided single-input single-output (SISO) communication system in [1], and proposed a beamforming design for RIS optimization that incorporates specified reradiation masks at the design stage.

*B. Main Contributions*

As previously stated, in our preliminary work [1], we introduced a beamforming design for RIS-aided SISO communication system that incorporates specified reradiation masks at the design stage using two methods: semidefinite programming (SDP) and a model-aided neural network architecture.

Although some prior studies, including our work in [1], have addressed reradiation mask constraints, they are limited to the SISO case, as discussed in detail in the comprehensive literature review presented in [1]. In this paper, we extend our study to an RIS-aided multi-stream MIMO communication system. We study the joint design of the transmit precoders matrices and the RIS phase shift vector, aiming to maximize the minimum achievable rate of multiple receivers, which introduces additional challenges in the optimization. To address the formulated non-convex max-min optimization problem, we first develop an alternating approach that converts the sub-problems into convex ones. To further enhance the efficiency of our approach, we develop a neural network architecture based on the objective function of the max-min optimization problem. Furthermore, practical RIS restrictions, e.g., discrete phase shifts, are considered. To elaborate, our contributions are summarized as follows.

- First, we formulate a max-min problem to maximize the minimum achievable rate among multiple receivers in the system, considering transmit power and reradiation mask constraints. We examine two case studies: *i)* unit amplitude and continuous phase (UACP) and *ii)* unit amplitude and discrete phase (UADP) for the RIS scattering coefficients, which are commonly used in the literature [50]. We use the Arimoto-Blahut structure to simplify the expression of the achievable rate. Based on this, we propose an alternating optimization approach to find a suboptimal solution for the transmit precoding matrices and the RIS phase shift vector in the UACP case by solving sub-problems constituted by quadratic programs with quadratic constraints (QPQC). The convergence and complexity of the proposed algorithm are also discussed.
- To improve the efficiency of the proposed approach, we develop a model-based neural network design. The neural network takes as inputs the angle of incidence and the desired angles of reflection, and outputs the concatenated vector of RIS phase shifts and the vectorized transmit precoder matrices. To efficiently represent the desired angles of incidence and reflection, we introduce a one-hot encoding method. Both the alternating-based optimization and the model-aided neural networks method account for specified reradiation constraints on the radiated power.
- We address practical RIS limitations by solving the optimization problem for discrete phase shifts using a heuristic approach based on a greedy search algorithm. Numerical simulations demonstrate that the proposed methods effectively shape the reradiation pattern of the RIS towards desired directions of reradiation while fulfilling specified reradiation mask constraints towards other directions of reradiation. The neural network-based method shows good performance with reduced computation time. The discrete phase shifts scheme performs well with only four phase shift levels, exhibiting a small reduction in terms of beamforming gain towards the desired directions of reradiation.

*Organization:* The remainder of the paper is organized as follows. Section II introduces the system model and Section III

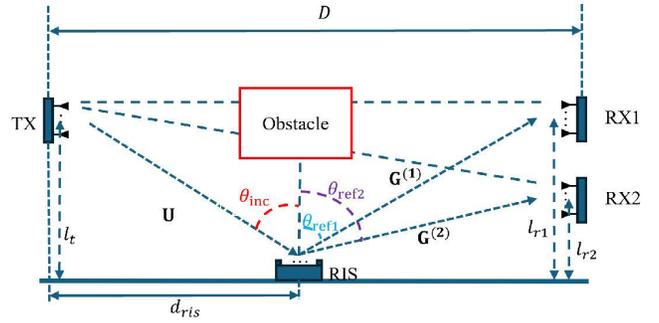

Fig. 1. Aerial view of the considered communication system.

formulates the optimization problem. An alternating method-based algorithm is proposed in Section IV, followed by the practical RIS case study with discrete phase shifts in Section V. The model-based neural network approach is addressed in Section VI. Numerical results are provided in Section VII to evaluate the performance of the proposed schemes. Finally, Section VIII concludes the paper.

*Notation:* The following notations are used throughout this paper. Matrices and vectors are denoted by bold capital and bold small fonts, respectively. $|z|$, $\mathcal{R}(z)$ and $\mathcal{J}(z)$ denote the absolute value, real and imaginary parts of $z$, respectively. $\mathbb{E}\{\cdot\}$ denotes the expectation. $\mathbf{I}_N$ denotes the identity matrix of size $N \times N$. $\mathbf{A}^{-1}$, $\mathbf{A}^T$, $\mathbf{A}^H$, $Tr(\mathbf{A})$ and $\det(\mathbf{A})$ denote the inverse, transpose, hermitian, trace and determinant of matrix $\mathbf{A}$, respectively. $\mathrm{diag}(\mathbf{a})$ denotes the diagonal matrix with the vector $\mathbf{a}$ on its diagonal. $\mathbf{A} \odot \mathbf{B}$ denotes the Hadamard product of $\mathbf{A}$ and $\mathbf{B}$. $\|\mathbf{a}\|$ denotes the $l_2$-norm of vector $\mathbf{a}$. $\cup$ is the union symbol. $\mathcal{CN}(\boldsymbol{\mu}, \boldsymbol{\Sigma})$ denotes the distribution of a circularly symmetric complex Gaussian (CSCG) random vector with mean vector $\boldsymbol{\mu}$ and covariance matrix $\boldsymbol{\Sigma}$; and $\sim$ stands for "distributed as". Finally, $\mathcal{O}(\cdot)$ stands for the big-O notation.

## II. SYSTEM MODEL

Consider an RIS-aided downlink communication system channel, as shown in Figure 1. The system includes a transmitter with $N_t$ antennas and two receivers, each equipped with $N_{r_1}$ and $N_{r_2}$ antennas, respectively. The proposed approach can be readily generalized to the multi-receiver and multi-beam case, as reported in the appendix for ease of exposition [2]. The distance between the transmitter and receivers is denoted by $D$. It is assumed that no direct link exists between the transmitter and receivers due to the presence of blocking objects. Consequently, communication occurs through an RIS with $N_{\mathrm{ris}}$ reflecting elements, arranged in a uniform rectangular array. The transmitter aims to transmit data to both users, hence the RIS must be optimized to generate two beams directed at the two intended users.

The midpoint of the RIS is at a distance $d_{ris}$ from the plane containing the transmitter, and the distance between the midpoint of the transmitter (respectively, receivers) and the plane containing the RIS is $l_t$ (respectively, $l_{r1}$ and $l_{r2}$). The separation between the centers of two adjacent RIS elements in both dimensions is $s_{\text{ris}} = \lambda/2$, where $\lambda$ is the wavelength. Also, $\theta_{inc}$ is the angle of incidence from the transmitter to the RIS, $\theta_{ref1}$ and $\theta_{ref2}$ are the angles of reflection from the RIS to receiver $r_1$ and receiver $r_2$, respectively [3]. The channel matrix between the transmitter and the RIS is denoted by $\mathbf{U} \in \mathbb{C}^{N_{\text{ris}} \times N_t}$, the channel matrix between the RIS and the receiver $r_1$ is denoted by $\mathbf{G}^{(1)} \in \mathbb{C}^{N_{r_1} \times N_{\text{ris}}}$, and the channel matrix between the RIS and the receiver $r_2$ is denoted by $\mathbf{G}^{(2)} \in \mathbb{C}^{N_{r_2} \times N_{\text{ris}}}$.

---
[2]The appendix provides the main proof for the proposed optimization method applied to the UACP case. The UADP case and the model-driven neural network optimization are straightforward extensions.

For the receiver $r_i$, $i = 1, 2$, the received signal is given by

$$\mathbf{y}_i = \mathbf{H}_i \mathbf{F}_i \mathbf{s}_i + \mathbf{H}_i \mathbf{F}_{\bar{i}} \mathbf{s}_{\bar{i}} + \mathbf{n}_i \quad (1)$$

where $\bar{i} = 3-i$, $i = 1, 2$. $\mathbf{H}_i \in \mathbb{C}^{N_{r_i} \times N_t}$ is the channel gain of the end-to-end transmitter-RIS-receiver $r_i$ link, $\mathbf{F}_i \in \mathbb{C}^{N_t \times N_{r_i}}$ is the transmit precoder matrix to the receiver $r_i$, $i = 1, 2$, and $\mathbf{s}_i \in \mathbb{C}^{N_{r_i} \times 1}$ is the transmit symbol vector with covariance matrix $\mathbf{I}_{N_{r_i}}$, assuming that one data stream is transmitted per receive antenna. The noise $\mathbf{n}_i \in \mathbb{C}^{N_{r_i} \times 1}$ consists of independent and identically distributed (i.i.d.) elements with distribution $\mathcal{CN}(0, \sigma^2 \mathbf{I}_{N_{r_i}})$, where $\sigma^2$ is the noise variance.

Considering the presence of the RIS, the channel $\mathbf{H}_i$, $i = 1, 2$ can be expressed as

$$\mathbf{H}_i = (\beta_{\text{INDIR}}^{(i)})^{-1/2} \mathbf{G}^{(i)} \text{diag}(\boldsymbol{\theta}) \mathbf{U} \quad (2)$$

where $\boldsymbol{\theta} = [\theta_1, \theta_2, \ldots, \theta_{N_{\text{ris}}}]^T \in \mathbb{C}^{N_{\text{ris}} \times 1}$ is the vector of RIS reflection coefficients, and $\text{diag}(\boldsymbol{\theta})$ returns a square matrix whose main diagonal contains the entries of $\boldsymbol{\theta}$. Also, $1/\beta_{\text{INDIR}}^{(i)}$ is the free-space path-loss factor, which is defined next.

We consider a Rician fading channel model [4]. Thus, the channel $\mathbf{U}$ between the transmitter and the RIS is given by

$$\mathbf{U} = \frac{1}{\sqrt{K_r + 1}} (\sqrt{K_r} \mathbf{U}_{\text{LOS}} + \mathbf{U}_{\text{NLOS}}) \quad (3)$$

where $\mathbf{U}_{\text{LOS}}(l, t) = e^{-j2\pi d_{l,t}/\lambda}$, and $d_{l,t}$ is the distance between the $t$-th transmit antenna and the $l$-th RIS element. The elements of $\mathbf{U}_{\text{NLOS}}$ are i.i.d. according to $\mathcal{CN}(0, 1)$. The Rician factor $K_r$ is chosen from the interval $[0, +\infty)$.

---

[3] In this work, the angles of incidence and reflection are assumed to be known. Estimating these angles in practice, particularly in the presence of multipath propagation, is a well-studied problem with several effective methods available in the literature [51]–[54]. However, integrating angle estimation into our framework is beyond the scope of this work. Numerical results in the presence of imperfect channel estimation are presented in Section VII.

[4] The Rician model is used since it reflects both the dominant line-of-sight (LOS) path and the scattered multipath from surrounding obstacles, typically present in wireless scenarios [9]–[12], [17]–[19]. In our model, we assume that the direct LOS path between the BS and the receivers is blocked, with the RIS providing the main propagation path. While weak NLOS components independent of the RIS may exist in practice, they are generally neglected due to their low power and limited impact. These NLOS paths are captured statistically through the flat Rician fading model, which is applied to the BS–RIS and RIS–receiver links. A more general model that includes both direct and RIS-aided links could also be considered and would remain compatible with the proposed optimization methods. Finally, we note that the proposed optimization approaches are independent of the specific channel model used.

Similarly, $\mathbf{G}^{(i)}$, $i = 1, 2$ can be expressed as

$$\mathbf{G}^{(i)} = \frac{1}{\sqrt{K_r + 1}} (\sqrt{K_r} \mathbf{G}_{\text{LOS}}^{(i)} + \mathbf{G}_{\text{NLOS}}^{(i)}) \quad (4)$$

where $\mathbf{G}_{\text{LOS}}^{(i)}(r, l) = e^{-j2\pi d_{r,l}/\lambda}$ and $d_{r,l}$ is the distance between the $l$-th RIS element and the $r$-th antenna of receiver $r_i$, $i = 1, 2$. The elements of $\mathbf{G}_{\text{NLOS}}^{(i)}$ are i.i.d. according to $\mathcal{CN}(0, 1)$.

In addition, we assume that the transmitter and the two receivers are in the far field of the RIS and of each other. Therefore, the free-space path-loss factor of the end-to-end link can be computed according to [55, Eq. (8)]

$$\frac{1}{\beta_{\text{INDIR}}^{(i)}} = \frac{G_t G_r \lambda^4}{256\pi^2} \frac{\cos \gamma_1 \cos \gamma_2^{(i)}}{d_1^2 (d_2^{(i)})^2} \quad (5)$$

where $d_1 = \sqrt{d_{\text{ris}}^2 + l_t^2}$ is the distance between the midpoint of the transmitter and the center-point of RIS, $d_2^{(i)} = \sqrt{(D - d_{\text{ris}})^2 + l_{r_i}^2}$ is the distance between the center-point of RIS and the midpoint of the receiver $r_i$, and $G_t = G_r = 2$ for simplicity. Also, $\gamma_1$ is the angle between the direction of the incident wave from the transmitter to the center-point of RIS and the vector normal to the RIS, and $\gamma_2^{(i)}$ is the angle between the vector normal to the RIS and the direction of the reflected wave from the center-point of RIS to the receiver. Hence, $\cos \gamma_1 = l_t/d_1$ and $\cos \gamma_2^{(i)} = l_{r_i}/d_2^{(i)}$, $i = 1, 2$.

## III. PROBLEM FORMULATION

We consider a max-min optimization problem which aims to maximize the minimum achievable rate among the two receivers $r_1$ and $r_2$, by jointly optimizing the transmit precoding matrices $\mathbf{F}_1$ and $\mathbf{F}_2$, and the phase shift vector $\boldsymbol{\theta}$. The problem can be formulated as

$$\max_{\boldsymbol{\theta}, \mathbf{F}_1, \mathbf{F}_2} \min\{R_1(\boldsymbol{\theta}, \mathbf{F}_1, \mathbf{F}_2), R_2(\boldsymbol{\theta}, \mathbf{F}_1, \mathbf{F}_2)\} \quad (6a)$$

$$\text{s.t.} \quad \theta_n \in \mathcal{F} \quad \forall n \in \{1, \cdots, N_{ris}\} \quad (6b)$$

$$\sum_{i=1}^{2} \text{Tr}(\mathbf{F}_i \mathbf{F}_i^H) \leq P_{max} \quad (6c)$$

$$\Pr(\boldsymbol{\theta}, \mathbf{F}_1, \mathbf{F}_2, \theta^{(ob)}) \leq \rho \quad \theta^{(ob)} \in \mathcal{A} \quad (6d)$$

where $R_i$, $i = 1, 2$ is the achievable rate of receiver $r_i$, which is defined as

$$R_i = \log_2 \left( \det \left( \mathbf{I}_{N_{r_i}} + \mathbf{H}_i \mathbf{F}_i \boldsymbol{\Omega}_{\bar{i}}^{-1} \mathbf{F}_i^H \mathbf{H}_i^H \right) \right) \quad (7)$$

where $\boldsymbol{\Omega}_{\bar{i}} = \mathbf{H}_i \mathbf{F}_{\bar{i}} \mathbf{F}_{\bar{i}}^H \mathbf{H}_i^H + \sigma^2 \mathbf{I}_{N_{r_i}}$, and $\bar{i} \triangleq 3 - i$, represents the interference-plus-noise covariance matrix for receiver $r_i$. $\mathcal{F}$ is the set of allowed phase shifts of the RIS.

To maximize the received power of both receivers, we consider a maximum signal reflection where the amplitude is equal to one, i.e., $|\theta_n| = 1, \forall n$. In what follows, we consider two cases for the feasibility set of $\theta_n$

i) *Unit Amplitude and Continuous Phase (UACP):* each phase shift $\theta_n$ can be adjusted to any desired continuous value, hence

$$\mathcal{F}_1 = \{\theta_n \,|\, \theta_n = e^{j\varphi_n}, \varphi_n \in [0, 2\pi)\} \quad (8)$$

ii) *Unit Amplitude and Discrete Evenly-spaced Phase (UADP):* due to practical limitations of the RIS, the reflecting elements can only have finite phase shift levels. We assume $L$ levels that are equally spaced within $[0, 2\pi)$, such as

$$\mathcal{F}_2 = \left\{ \theta_n \Big| \theta_n = e^{j\varphi_n}, \varphi_n \in \left\{ 0, \frac{2\pi}{L}, \ldots, \frac{2\pi(L-1)}{L} \right\} \right\} \quad (9)$$

The constraint (6c) represents the transmit power limit, where $P_{max}$ is the total transmit power for receiver $r_i$, $i = 1, 2$, and (6d) is the reradiation mask constraint where $\Pr(\boldsymbol{\theta}, \mathbf{F}_1, \mathbf{F}_2, \theta^{(ob)})$ is the power scattered by the RIS towards the direction $\theta^{(ob)}$, when the vector of reflection coefficients is set to $\boldsymbol{\theta}$, which is defined as

$$\begin{aligned}
&\Pr(\boldsymbol{\theta}, \mathbf{F}_1, \mathbf{F}_2, \theta^{(ob)}) \\
&= ||\mathbf{G}^{(ob)} \mathrm{diag}(\boldsymbol{\theta}) \mathbf{U} \mathbf{F}_1 \mathbf{s}_1 + \mathbf{G}^{(ob)} \mathrm{diag}(\boldsymbol{\theta}) \mathbf{U} \mathbf{F}_2 \mathbf{s}_2||^2 \\
&= \mathrm{Tr}\left( \mathbf{G}^{(ob)} \mathrm{diag}(\boldsymbol{\theta}) \mathbf{U} \mathbf{F}_1 \mathbf{F}_1^H \mathbf{U}^H \mathrm{diag}(\boldsymbol{\theta})^H (\mathbf{G}^{(ob)})^H \right) \\
&\quad + \mathrm{Tr}\left( \mathbf{G}^{(ob)} \mathrm{diag}(\boldsymbol{\theta}) \mathbf{U} \mathbf{F}_2 \mathbf{F}_2^H \mathbf{U}^H \mathrm{diag}(\boldsymbol{\theta})^H (\mathbf{G}^{(ob)})^H \right)
\end{aligned} \quad (10)$$

where $\theta^{(ob)}$ denotes the generic direction towards which the reradiated power must be below the threshold $\rho$ [5], $\mathbf{G}^{(ob)}$ is the corresponding channel matrix from the RIS to the generic observation angle.

Unfortunately, the optimization problem (6) is challenging to solve due to its non-convex nature. Following the approach in [23], we utilize the Arimoto-Blahut structure to reformulate the expression for the achievable rate as

$$R_i = \max_{q(\mathbf{s}_i|\mathbf{y}_i)} \mathbb{E} \left[ \log_2 \left( \frac{q(\mathbf{s}_i \mid \mathbf{y}_i)}{\mathcal{CN}(\mathbf{0}, \mathbf{I}_{N_{r_i}})} \right) \right] \quad (11)$$

According to [23], the optimal $q^\star(\mathbf{s}_i \mid \mathbf{y}_i)$ is the posterior probability $p(\mathbf{s}_i \mid \mathbf{y}_i)$, and according to [56, Theorem 10.3, p. 326], $p(\mathbf{s}_i \mid \mathbf{y}_i)$ can be derived such that it follows the complex Gaussian distribution $\mathcal{CN}(\mathbf{W}_i^\star \mathbf{y}_i, \boldsymbol{\Sigma}_i^\star)$ with

$$\mathbf{W}_i^\star = (\mathbf{H}_i \mathbf{F}_i)^H \left( (\mathbf{H}_i \mathbf{F}_i)(\mathbf{H}_i \mathbf{F}_i)^H + \boldsymbol{\Omega}_{\bar{i}} \right)^{-1} \quad (12)$$

$$\boldsymbol{\Sigma}_i^\star = \mathbf{I}_{N_{r_i}} - \mathbf{W}_i^\star \mathbf{H}_i \mathbf{F}_i \quad (13)$$

where $\boldsymbol{\Omega}_{\bar{i}} = \mathbf{G}^{(i)} \mathrm{diag}(\boldsymbol{\theta}) \mathbf{U} \mathbf{F}_{\bar{i}} \mathbf{F}_{\bar{i}}^H \mathbf{U}^H \mathrm{diag}(\boldsymbol{\theta})^H (\mathbf{G}^{(i)})^H + \sigma^2 \mathbf{I}_{N_{r_i}}$. According to (11)-(13), the optimization problem in (6) can be reformulated as shown in (14) at the bottom of the next page.

---

[5]The reradiation mask constraint is applied at the RIS, which lacks active components and signal processing capabilities and is more susceptible to power leakage in unintended directions [3], [45]. On the other hand, the BS can use techniques such as Dolph-Chebyshev weighting or amplitude tapering to suppress the sidelobes. Additionally, since the BS–receiver links are blocked in our model, the intended signal from the BS is not part of the analysis. A possible generalization of this work could include radiation mask constraints at the BS. However, this is not considered here for simplicity, as the BS already has more advanced beam shaping capabilities compared to the RIS. Jointly optimizing the BS beam pattern along with the RIS configuration could be a promising direction for future work.

By computing the expectation term in (14a), the objective function can be written as

$$\begin{aligned}
&\mathbb{E} \left[ \log_2 \left( \frac{\mathcal{CN}(\mathbf{W}_i \mathbf{y}_i, \boldsymbol{\Sigma}_i)}{\mathcal{CN}(\mathbf{0}, \mathbf{I}_{N_{r_i}})} \right) \right] \\
&= 2 \, \mathrm{Re} \left( \mathrm{Tr} \left( \boldsymbol{\Sigma}_i^{-1} \mathbf{W}_i \mathbf{G}^{(i)} \mathrm{diag}(\boldsymbol{\theta}) \mathbf{U} \mathbf{F}_i \right) \right) - \mathrm{Tr}(\boldsymbol{\Sigma}_i^{-1}) \\
&\quad - \mathrm{Tr} \Big( \mathbf{F}_i^H \mathbf{U}^H \mathrm{diag}(\boldsymbol{\theta})^H \left( \mathbf{G}^{(i)} \right)^H \mathbf{W}_i^H \boldsymbol{\Sigma}_i^{-1} \\
&\hspace{6em} \mathbf{W}_i \mathbf{G}^{(i)} \mathrm{diag}(\boldsymbol{\theta}) \mathbf{U} \mathbf{F}_i \Big) \\
&\quad - \mathrm{Tr} \Big( \mathbf{F}_{\bar{i}}^H \mathbf{U}^H \mathrm{diag}(\boldsymbol{\theta})^H \left( \mathbf{G}^{(i)} \right)^H \mathbf{W}_i^H \boldsymbol{\Sigma}_i^{-1} \\
&\hspace{6em} \mathbf{W}_i \mathbf{G}^{(i)} \mathrm{diag}(\boldsymbol{\theta}) \mathbf{U} \mathbf{F}_{\bar{i}} \Big) \\
&\quad - \sigma^2 \mathrm{Tr}\left( \mathbf{W}_i^H \boldsymbol{\Sigma}_i^{-1} \mathbf{W}_i \right) - N_{r_i} \log_2(\det(\boldsymbol{\Sigma}_i)) + N_{r_i}
\end{aligned} \quad (15)$$

In what follows, we address the optimization problem in (14) using the AO method. This involves iteratively solving sub-problems with respect to one variable at a time while keeping the other variables fixed.

## IV. Optimization Algorithm for UACP

In this section, we initially consider that the phase shift of each RIS element is continuous, with each element having an amplitude equal to one, i.e., $\theta_n \in \mathcal{F}_1$. We then solve the problem in (14) using the AO approach.

### A. Update $\mathbf{W}_1$, $\boldsymbol{\Sigma}_1$, $\mathbf{W}_2$, and $\boldsymbol{\Sigma}_2$

Here, we optimize $\mathbf{W}_i$ and $\boldsymbol{\Sigma}_i$ while keeping $\mathbf{F}_i$ and $\boldsymbol{\theta}$ fixed. The solutions for these matrices are provided in (12) and (13) of [23]. Therefore, they are not reported for brevity.

### B. Update the RIS Phase Shift Vector $\boldsymbol{\theta}$

Here, we optimize $\boldsymbol{\theta}$ keeping $\mathbf{F}_i$, $\mathbf{W}_i$ and $\boldsymbol{\Sigma}_i$, $i = 1, 2$, fixed. To this end, we note that the first, third, and fourth terms of the achievable rate in (15) depend on the phase shift vector $\boldsymbol{\theta}$.

The first term in (15) can be written as

$$\begin{aligned}
&2 \, \mathrm{Re} \left( \mathrm{Tr} \left( \boldsymbol{\Sigma}_i^{-1} \mathbf{W}_i \mathbf{G}^{(i)} \mathrm{diag}(\boldsymbol{\theta}) \mathbf{U} \mathbf{F}_i \right) \right) \\
&= 2 \, \mathrm{Re} \left( \mathrm{Tr} \left( \mathbf{U} \mathbf{F}_i \boldsymbol{\Sigma}_i^{-1} \mathbf{W}_i \mathbf{G}^{(i)} \mathrm{diag}(\boldsymbol{\theta}) \right) \right) \\
&= 2 \, \mathrm{Re} \left( \left( \boldsymbol{\theta}^H \mathbf{b}_i \right) \right)
\end{aligned} \quad (16)$$

where $\mathbf{b}_i = \mathrm{diag}\left( \left( \mathbf{U} \mathbf{F}_i \boldsymbol{\Sigma}_i^{-1} \mathbf{W}_i \mathbf{G}^{(i)} \right)^H \right)$.

Let us define $\mathbf{A}_i = \left( \mathbf{G}^{(i)} \right)^H \mathbf{W}_i^H \boldsymbol{\Sigma}_i^{-1} \mathbf{W}_i \mathbf{G}^{(i)}$, $\mathbf{B}_i = \left( \mathbf{U} \mathbf{F}_i \mathbf{F}_i^H \mathbf{U}^H \right)^T$, $\mathbf{C}_i = \left( \mathbf{G}^{(i)} \right)^H \mathbf{W}_i^H \boldsymbol{\Sigma}_i^{-1} \mathbf{W}_i \mathbf{G}^{(i)}$ and $\mathbf{D}_i = \left( \mathbf{U} \mathbf{F}_{\bar{i}} \mathbf{F}_{\bar{i}}^H \mathbf{U}^H \right)^T$. The third term in (15) can be written as

$$\begin{aligned}
&-\mathrm{Tr}\Big( \mathbf{F}_i^H \mathbf{U}^H \mathrm{diag}(\boldsymbol{\theta})^H (\mathbf{G}^{(i)})^H \mathbf{W}_i^H \boldsymbol{\Sigma}_i^{-1} \mathbf{W}_i \mathbf{G}^{(i)} \mathrm{diag}(\boldsymbol{\theta}) \mathbf{U} \mathbf{F}_i \Big) \\
&= -\mathrm{Tr}\Big( \mathrm{diag}(\boldsymbol{\theta})^H (\mathbf{G}^{(i)})^H \mathbf{W}_i^H \boldsymbol{\Sigma}_i^{-1} \mathbf{W}_i \mathbf{G}^{(i)} \mathrm{diag}(\boldsymbol{\theta}) \mathbf{U} \mathbf{F}_i \mathbf{F}_i^H \mathbf{U}^H \Big) \\
&= -\boldsymbol{\theta}^H (\mathbf{A}_i \odot \mathbf{B}_i) \boldsymbol{\theta}
\end{aligned} \quad (17)$$

where the last equality comes from $\mathrm{Tr}\left( \mathrm{diag}(\boldsymbol{\theta})^H \mathbf{A}_i \mathrm{diag}(\boldsymbol{\theta}) \mathbf{B}_i \right) = \boldsymbol{\theta}^H (\mathbf{A}_i \odot \mathbf{B}_i) \boldsymbol{\theta}$. Similarly, the fourth term in (15) can be expressed as

$$\begin{aligned}
&-\mathrm{Tr}\Big( \mathbf{F}_{\bar{i}}^H \mathbf{U}^H \mathrm{diag}(\boldsymbol{\theta})^H (\mathbf{G}^{(i)})^H \mathbf{W}_i^H \boldsymbol{\Sigma}_i^{-1} \mathbf{W}_i \mathbf{G}^{(i)} \mathrm{diag}(\boldsymbol{\theta}) \mathbf{U} \mathbf{F}_{\bar{i}} \Big) \\
&= -\mathrm{Tr}\Big( \mathrm{diag}(\boldsymbol{\theta})^H (\mathbf{G}^{(i)})^H \mathbf{W}_i^H \boldsymbol{\Sigma}_i^{-1} \mathbf{W}_i \mathbf{G}^{(i)} \mathrm{diag}(\boldsymbol{\theta}) \mathbf{U} \mathbf{F}_{\bar{i}} \mathbf{F}_{\bar{i}}^H \mathbf{U}^H \Big) \\
&= -\boldsymbol{\theta}^H (\mathbf{C}_i \odot \mathbf{D}_i) \boldsymbol{\theta}
\end{aligned} \quad (18)$$

The remaining terms in (15) that are not related to $\boldsymbol{\theta}$ can be considered as constant, and are denoted by

$$c_i = -\sigma^2 \text{Tr}(\mathbf{W}_i^H \boldsymbol{\Sigma}_i^{-1} \mathbf{W}_i) - N_{r_i} \log_2(\det(\boldsymbol{\Sigma}_i)) + N_{r_i} - \text{Tr}(\boldsymbol{\Sigma}_i^{-1}) \quad (19)$$

The mask constraint in (14d) can also be simplified by using similar algebraic manipulations, as follows

$$\begin{aligned}
\Pr&(\boldsymbol{\theta}, \mathbf{F}_1, \mathbf{F}_2, \theta^{(ob)}) \\
&= \text{Tr}\left(\mathbf{G}^{(ob)} \text{diag}(\boldsymbol{\theta}) \mathbf{U} \mathbf{F}_1 \mathbf{F}_1^H \mathbf{U}^H \text{diag}(\boldsymbol{\theta})^H (\mathbf{G}^{(ob)})^H\right) \\
&\quad + \text{Tr}\left(\mathbf{G}^{(ob)} \text{diag}(\boldsymbol{\theta}) \mathbf{U} \mathbf{F}_2 \mathbf{F}_2^H \mathbf{U}^H \text{diag}(\boldsymbol{\theta})^H (\mathbf{G}^{(ob)})^H\right) \\
&= \text{Tr}\left(\text{diag}(\boldsymbol{\theta})^H (\mathbf{G}^{(ob)})^H \mathbf{G}^{(ob)} \text{diag}(\boldsymbol{\theta}) \mathbf{U} \mathbf{F}_1 \mathbf{F}_1^H \mathbf{U}^H\right) \\
&\quad + \text{Tr}\left(\text{diag}(\boldsymbol{\theta})^H (\mathbf{G}^{(ob)})^H \mathbf{G}^{(ob)} \text{diag}(\boldsymbol{\theta}) \mathbf{U} \mathbf{F}_2 \mathbf{F}_2^H \mathbf{U}^H\right) \\
&= \boldsymbol{\theta}^H \left(\mathbf{Q}^{(ob)} \odot \mathbf{T}_1 + \mathbf{Q}^{(ob)} \odot \mathbf{T}_2\right) \boldsymbol{\theta} \quad (20)
\end{aligned}$$

where $\mathbf{Q}^{(ob)} = (\mathbf{G}^{(ob)})^H \mathbf{G}^{(ob)}$, $\mathbf{T}_1 = \left(\mathbf{U} \mathbf{F}_1 \mathbf{F}_1^H \mathbf{U}^H\right)^T$ and $\mathbf{T}_2 = \left(\mathbf{U} \mathbf{F}_2 \mathbf{F}_2^H \mathbf{U}^H\right)^T$.

The constraint in (14b), where $\theta_n \in \mathcal{F}_1$, can be written in a quadratic form as

$$\boldsymbol{\theta}^H \mathbf{I}_{N_{ris}}(:,n) \left(\mathbf{I}_{N_{ris}}(:,n)\right)^H \boldsymbol{\theta} = 1, \quad \forall n \in \{1, \ldots, N_{ris}\} \quad (21)$$

where $\mathbf{I}_{N_{ris}}(:,n)$ is an $N_{ris} \times 1$ vector with entries equal to one at the $n$-th position and zero elsewhere.

Substituting (16)-(18) into (14a), (20) and (21) in (14d), and (14b) respectively, and removing terms irrelevant to $\boldsymbol{\theta}$, the sub-problem for optimizing with respect to $\boldsymbol{\theta}$ can be reformulated as

$$\begin{aligned}
\max_{\boldsymbol{\theta}} \quad & \min\left\{-\boldsymbol{\theta}^H \mathbf{E}_1 \boldsymbol{\theta} + 2\,\text{Re}\left(\boldsymbol{\theta}^H \mathbf{b}_1\right) + c_1,\right. \\
& \quad \left.-\boldsymbol{\theta}^H \mathbf{E}_2 \boldsymbol{\theta} + 2\,\text{Re}\left(\boldsymbol{\theta}^H \mathbf{b}_2\right) + c_2\right\} \\
\text{s.t.} \quad & \boldsymbol{\theta}^H \mathbf{I}_{N_{ris}}(:,n) \mathbf{I}_{N_{ris}}(:,n)^H \boldsymbol{\theta} = 1 \quad \forall n \\
& \boldsymbol{\theta}^H (\mathbf{Q}^{(ob)} \odot \mathbf{T}_1 + \mathbf{Q}^{(ob)} \odot \mathbf{T}_2) \boldsymbol{\theta} \le \rho \quad \theta^{(ob)} \in \mathcal{A}
\end{aligned} \quad (22)$$

where $\mathbf{E}_1 = (\mathbf{A}_1 \odot \mathbf{B}_1 + \mathbf{C}_1 \odot \mathbf{D}_1)$ and $\mathbf{E}_2 = (\mathbf{A}_2 \odot \mathbf{B}_2 + \mathbf{C}_2 \odot \mathbf{D}_2)$.

It it not difficult to verify that $\mathbf{Q}^{(ob)}$, $\mathbf{A}_i$, $\mathbf{B}_i$, $\mathbf{C}_i$, $\mathbf{D}_i$, and $\mathbf{T}_i$ $i = 1, 2$, are Hermitian semi-positive definite matrices. Consequently, $\mathbf{E}_i$ and $\left(\mathbf{Q}^{(ob)} \odot \mathbf{T}_1\right)$ are also semi-positive definite matrices. Thus, the problem in (22) is a convex QCQP, which can be solved efficiently by using typical convex optimization solvers, such as CVX.

### C. Update the Transmit Precoding $\mathbf{F}_1$ and $\mathbf{F}_2$

Here, we optimize $\mathbf{F}_1$ and $\mathbf{F}_2$ while keeping $\boldsymbol{\theta}$, $\mathbf{W}_i$ and $\boldsymbol{\Sigma}_i$ fixed for $i = 1, 2$. In (15), the first, third and fourth terms depend on the transmit precoding matrices $\mathbf{F}_1$ and $\mathbf{F}_2$. Let us define $\mathbf{K}_i = \boldsymbol{\Sigma}_i^{-1} \mathbf{W}_i \mathbf{G}^{(i)} \text{diag}(\boldsymbol{\theta}) \mathbf{U}$ and $\mathbf{J}_i = \mathbf{U}^H \text{diag}(\boldsymbol{\theta})^H (\mathbf{G}^{(i)})^H \mathbf{W}_i^H \boldsymbol{\Sigma}_i^{-1} \mathbf{W}_i \mathbf{G}^{(i)} \text{diag}(\boldsymbol{\theta}) \mathbf{U}$.

For simplicity, we initially disregard the mask constraint (14d). Then, by updating $\boldsymbol{\theta}$ in the sub-problem (22), we ensure that the mask constraint is fulfilled. By ignoring the terms irrelevant to $\mathbf{F}_1$ and $\mathbf{F}_2$, the optimization sub-problem is formulated as

$$\begin{aligned}
\max_{\mathbf{F}_1, \mathbf{F}_2} \quad & \min\{M_1(\mathbf{F}_1, \mathbf{F}_2), M_2(\mathbf{F}_1, \mathbf{F}_2)\} \\
\text{s.t.} \quad & \sum_{i=1}^{2} \text{Tr}\left(\mathbf{F}_i \mathbf{F}_i^H\right) \le P_{max}
\end{aligned} \quad (23)$$

where $M_1(\mathbf{F}_1, \mathbf{F}_2) = -\text{Tr}\left(\mathbf{F}_1^H \mathbf{J}_1 \mathbf{F}_1\right) - \text{Tr}\left(\mathbf{F}_2^H \mathbf{J}_1 \mathbf{F}_2\right) + 2\text{Re}\left(\text{Tr}\left(\mathbf{F}_1^H \mathbf{K}_1^H\right)\right) + c_1$, and $M_2(\mathbf{F}_1, \mathbf{F}_2) = -\text{Tr}\left(\mathbf{F}_1^H \mathbf{J}_2 \mathbf{F}_1\right) - \text{Tr}\left(\mathbf{F}_2^H \mathbf{J}_2 \mathbf{F}_2\right) + 2\text{Re}\left(\text{Tr}\left(\mathbf{F}_2^H \mathbf{K}_2^H\right)\right) + c_2$.

The sub-problem in (23) is still non-convex and challenging to solve. To address it, we decompose it by solving for $\mathbf{F}_1$ and $\mathbf{F}_2$ separately. For example, focusing on $\mathbf{F}_1$, the optimization problem is reformulated as follows:

$$\begin{aligned}
\max_{\mathbf{F}_1} \quad & \min\{-\text{Tr}\left(\mathbf{F}_1^H \mathbf{J}_1 \mathbf{F}_1\right) + 2\text{Re}\left(\text{Tr}\left(\mathbf{F}_1^H \mathbf{K}_1^H\right)\right) - v_1, \\
& \quad -\text{Tr}\left(\mathbf{F}_1^H \mathbf{J}_2 \mathbf{F}_1\right) + o_1\} \\
\text{s.t.} \quad & \text{Tr}\left(\mathbf{F}_1 \mathbf{F}_1^H\right) \le P_{max}^{(1)}
\end{aligned} \quad (24)$$

where $v_1 = \text{Tr}\left(\mathbf{F}_2^H \mathbf{J}_1 \mathbf{F}_2\right) - c_1$, $o_1 = -\text{Tr}\left(\mathbf{F}_2^H \mathbf{J}_2 \mathbf{F}_2\right) + 2\text{Re}\left(\text{Tr}\left(\mathbf{F}_2^H \mathbf{K}_2^H\right)\right) + c_2$ and $P_{max}^{(1)} = P_{max} - \text{Tr}\left(\mathbf{F}_2 \mathbf{F}_2^H\right)$.

It is not difficult to verify that $\mathbf{J}_1$ and $\mathbf{J}_2$ are semi-positive definite matrices. Hence, (24) is a convex QCQP and can be solved using CVX, similarly to (22). However, since CVX cannot directly handle matrix variables, we vectorize $\mathbf{F}_1$ and transform other corresponding matrices accordingly as

$$\mathbf{f}_1 = [\mathbf{F}_1(:,1)^T, \cdots, \mathbf{F}_1(:,N_{r_1})^T]^T \quad (25)$$
$$\mathbf{k}_1 = [\mathbf{K}_1(1,:), \cdots, \mathbf{K}_1(N_{r_1},:)] \quad (26)$$
$$\tilde{\mathbf{J}}_i = \text{blkdiag}(\underbrace{\mathbf{J}_i, \cdots, \mathbf{J}_i}_{N_{r_1}}), \quad i = 1, 2 \quad (27)$$

where $\text{blkdiag}(\cdot)$ returns the block diagonal matrix created by aligning the input matrices, such as

$$\tilde{\mathbf{J}}_i = \begin{bmatrix} \mathbf{J}_i & 0 & 0 \\ 0 & \ddots & 0 \\ 0 & 0 & \mathbf{J}_i \end{bmatrix} \quad (28)$$

---

$$\max_{\boldsymbol{\theta}, \mathbf{F}_1, \mathbf{F}_2, \mathbf{W}_1, \mathbf{W}_2, \boldsymbol{\Sigma}_1, \boldsymbol{\Sigma}_2} \min\left\{\mathbb{E}\left[\log_2\left(\frac{\mathcal{CN}(\mathbf{W}_1 \mathbf{y}_1, \boldsymbol{\Sigma}_1)}{\mathcal{CN}(\mathbf{0}, \mathbf{I}_{N_{r_1}})}\right)\right], \mathbb{E}\left[\log_2\left(\frac{\mathcal{CN}(\mathbf{W}_2 \mathbf{y}_2, \boldsymbol{\Sigma}_2)}{\mathcal{CN}(\mathbf{0}, \mathbf{I}_{N_{r_2}})}\right)\right]\right\} \quad (14a)$$

$$\text{s.t.} \quad \theta_n \in \mathcal{F} \quad \forall n \in \{1, \cdots, N_{ris}\} \quad (14b)$$

$$\sum_{i=1}^{2} \text{Tr}(\mathbf{F}_i \mathbf{F}_i^H) \le P_{max} \quad (14c)$$

$$\Pr(\boldsymbol{\theta}, \mathbf{F}_1, \mathbf{F}_2, \theta^{(ob)}) \le \rho \quad \theta^{(ob)} \in \mathcal{A} \quad (14d)$$

The problem in (24) can be re-expressed as

$$\max_{\mathbf{f}_1} \min\{-\mathrm{Tr}\left(\mathbf{f}_1^H \tilde{\mathbf{J}}_1 \mathbf{f}_1\right) + 2\mathrm{Re}\left(\mathrm{Tr}\left(\mathbf{f}_1^H \mathbf{k}_1^H\right)\right) - v_1,$$
$$-\mathrm{Tr}\left(\mathbf{f}_1^H \tilde{\mathbf{J}}_2 \mathbf{f}_1\right) + o_1\}$$
$$\text{s.t.} \quad \mathrm{Tr}\left(\mathbf{f}_1^H \mathbf{f}_1\right) \leq P_{max}^{(1)} \quad (29)$$

Hence, the sub-problem becomes a convex QCQP, which can be solved using CVX. Once the optimal solution is obtained, it is necessary to convert the optimal precoding vector $\mathbf{f}_1^\star \in \mathbb{C}^{N_t N_{r_1} \times 1}$ into a matrix form, i.e., $\mathbf{F}_1^\star \in \mathbb{C}^{N_t \times N_{r_1}}$ by inverting (25) as

$$\mathbf{F}_1^\star = [\mathbf{f}_1^\star(1:N_t), \cdots, \mathbf{f}_1^\star((N_{r_1}-1)N_t + 1 : N_t N_{r_1})] \quad (30)$$

The optimal precoding matrix $\mathbf{F}_2^\star$ can be obtained in a similar manner. The corresponding sub-problem is expressed as

$$\max_{\mathbf{F}_2} \min\{-\mathrm{Tr}\left(\mathbf{F}_2^H \mathbf{J}_1 \mathbf{F}_2\right) + o_2,$$
$$-\mathrm{Tr}\left(\mathbf{F}_2^H \mathbf{J}_2 \mathbf{F}_2\right) + 2\mathrm{Re}\left(\mathrm{Tr}\left(\mathbf{F}_2^H \mathbf{K}_2^H\right)\right) - v_2\}$$
$$\text{s.t.} \quad \mathrm{Tr}\left(\mathbf{F}_2 \mathbf{F}_2^H\right) \leq P_{max}^{(2)} \quad (31)$$

where $v_2 = \mathrm{Tr}\left(\mathbf{F}_1^H \mathbf{J}_2 \mathbf{F}_1\right) - c_2$, $o_2 = -\mathrm{Tr}\left(\mathbf{F}_1^H \mathbf{J}_1 \mathbf{F}_1\right) + 2\mathrm{Re}\left(\mathrm{Tr}\left(\mathbf{F}_1^H \mathbf{K}_1^H\right)\right) + c_1$ and $P_{max}^{(2)} = P_{max} - \mathrm{Tr}\left(\mathbf{F}_1 \mathbf{F}_1^H\right)$.

Next, we vectorize the matrix $\mathbf{F}_2$ and transform the corresponding matrices as

$$\mathbf{f}_2 = [\mathbf{F}_2(:,1)^T, \cdots, \mathbf{F}_2(:,N_{r_2})^T]^T \quad (32)$$
$$\mathbf{k}_2 = [\mathbf{K}_2(1,:), \cdots, \mathbf{K}_2(N_{r_2},:)] \quad (33)$$

Then, (31) can be rewritten as

$$\max_{\mathbf{f}_2} \min\{-\mathrm{Tr}\left(\mathbf{f}_2^H \tilde{\mathbf{J}}_1 \mathbf{f}_2\right) + o_2,$$
$$-\mathrm{Tr}\left(\mathbf{f}_2^H \tilde{\mathbf{J}}_2 \mathbf{f}_2\right) + 2\mathrm{Re}\left(\mathrm{Tr}\left(\mathbf{f}_2^H \mathbf{k}_2^H\right)\right) - v_2\}$$
$$\text{s.t.} \quad \mathrm{Tr}\left(\mathbf{f}_2^H \mathbf{f}_2\right) \leq P_{max}^{(2)} \quad (34)$$

which is a convex QCQP that can be solved using CVX. The optimal precoding vector $\mathbf{f}_2^\star \in \mathbb{C}^{N_t N_{r_2} \times 1}$ is then reshaped into its matrix form as follows

$$\mathbf{F}_2^\star = [\mathbf{f}_2^\star(1:N_t), \cdots, \mathbf{f}_2^\star((N_{r_2}-1)N_t + 1 : N_t N_{r_2})] \quad (35)$$

The complete AO algorithm used to solve (6) is summarized in Algorithm 1, where $f$ denotes the objective function in (14a).

### D. Convergence and Complexity

We analyze the convergence of the proposed AO algorithm where the original problem in (6) is solved iteratively. In each iteration, the original problem is decomposed into the three sub-problems in (22), (29) and (34). All three sub-problems are convex, ensuring the convergence of each of them individually. This iterative process leads to the following sequence of inequalities

$$f\left(\boldsymbol{\theta}^{(t)}, \mathbf{F}_1^{(t)}, \mathbf{F}_2^{(t)}\right) \overset{a}{\leq} f\left(\boldsymbol{\theta}^{(t)}, \mathbf{F}_1^{(t+1)}, \mathbf{F}_2^{(t)}\right) \quad (36)$$
$$\overset{b}{\leq} f\left(\boldsymbol{\theta}^{(t)}, \mathbf{F}_1^{(t+1)}, \mathbf{F}_2^{(t+1)}\right)$$
$$\overset{c}{\leq} f\left(\boldsymbol{\theta}^{(t+1)}, \mathbf{F}_1^{(t+1)}, \mathbf{F}_2^{(t+1)}\right)$$

**Algorithm 1** AO algorithm for solving problem (6)

1: **Input:** $t=0$, $\epsilon > 0$, $\boldsymbol{\theta}^{(0)}$, $\mathbf{F}_1^{(0)}$, $\mathbf{F}_2^{(0)}$;
2: Compute $f^{(0)} = f(\boldsymbol{\theta}^{(0)}, \mathbf{F}_1^{(0)}, \mathbf{F}_2^{(0)})$ from (14a);
3: **repeat**
4:     Compute $\mathbf{W}_1, \mathbf{W}_2$ from (12) and $\boldsymbol{\Sigma}_1, \boldsymbol{\Sigma}_2$ from (13);
5:     Compute $\mathbf{F}_1^{(t+1)}$ by solving (29);
6:     Update $\mathbf{W}_1, \mathbf{W}_2$ and $\boldsymbol{\Sigma}_1, \boldsymbol{\Sigma}_2$ with $\mathbf{F}_1^{(t+1)}$ in (12) and (13) respectively;
7:     Compute $\mathbf{F}_2^{(t+1)}$ by solving (34);
8:     Update $\mathbf{W}_1, \mathbf{W}_2$ and $\boldsymbol{\Sigma}_1, \boldsymbol{\Sigma}_2$ with $\mathbf{F}_2^{(t+1)}$ in (12) and (13) respectively;
9:     Compute $\boldsymbol{\theta}^{(t+1)}$ by solving (22);
10:    Update $f^{(t+1)} \leftarrow f\left(\boldsymbol{\theta}^{(t+1)}, \mathbf{F}_1^{(t+1)}, \mathbf{F}_2^{(t+1)}\right)$;
11:    **if** $f^{(t+1)} < f^{(t)}$ **then**
12:       $\boldsymbol{\theta}^{(t+1)} \leftarrow \boldsymbol{\theta}^{(t)}$, $\mathbf{F}_1^{(t+1)} \leftarrow \mathbf{F}_1^{(t)}$, $\mathbf{F}_2^{(t+1)} \leftarrow \mathbf{F}_2^{(t)}$;
13:    $t \leftarrow t+1$;
14: **until** $|f^{(t)} - f^{(t-1)}| \leq \epsilon$.

where the inequalities (a)-(c) hold since solving each sub-problem in (22), (29) and (34) leads to an improvement over the previous solution. As a result, the values of (36) form a monotonically non-decreasing sequence throughout the iterations [57]. Moreover, the function $f(\boldsymbol{\theta}, \mathbf{F}_1, \mathbf{F}_2)$ is upper-bounded by a finite value due to the constraint on the BS transmit power in (6c) and the finite number and nearly-passive assumption for the reflecting elements of the RIS, which ensures that the sequence converges to a finite value. Thus, the proposed AO algorithm is guaranteed to converge.

Next, we present the computational complexity. The updates of $\mathbf{W}_i$ and $\boldsymbol{\Sigma}_i$ are matrix operations. The computational complexity of obtaining the matrix operation $\mathbf{W}_i = (\mathbf{H}_i \mathbf{F}_i)^H \left((\mathbf{H}_i \mathbf{F}_i)(\mathbf{H}_i \mathbf{F}_i)^H + \boldsymbol{\Omega}_i\right)^{-1}$ is obtained as follows. Multiplying $\mathbf{H}_i \in \mathbb{C}^{N_{r_i} \times N_t}$ with $\mathbf{F}_i \in \mathbb{C}^{N_t \times N_{r_i}}$ has complexity $\mathcal{O}(N_{r_i}^2 N_t)$. The most computationally expensive operations are the matrix inversion and multiplication, and each of them has complexity $\mathcal{O}(N_{r_i}^3)$. Thus, the total complexity is $\mathcal{O}(N_{r_i}^2 N_t + N_{r_i}^3)$. The expression $\boldsymbol{\Sigma}_i = \left(\mathbf{I}_{N_{r_i}} - \mathbf{W}_i \mathbf{H}_i \mathbf{F}_i\right)$ involves subtracting the identity matrix from the product of three matrices. The complexity of the matrix multiplication $\mathbf{W}_i \mathbf{H}_i \mathbf{F}_i$ is $\mathcal{O}(N_{r_i}^2 N_t)$, given that $\mathbf{W}_i \in \mathbb{C}^{N_{r_i} \times N_{r_i}}$. The matrix subtraction $\left(\mathbf{I}_{N_{r_i}} - \mathbf{W}_i \mathbf{H}_i \mathbf{F}_i\right)$ has complexity $\mathcal{O}(N_{r_i}^2)$. Thus, the total complexity is $\mathcal{O}(N_{r_i}^2 N_t)$, and the overall complexity to update $\mathbf{W}_1, \mathbf{W}_2, \boldsymbol{\Sigma}_1$ and $\boldsymbol{\Sigma}_2$ is $\mathcal{O}(N_{r_1}^2 N_t + N_{r_1}^3 + N_{r_2}^2 N_t + N_{r_2}^3)$. Each QCQP sub-problem in (22), (29) and (34) is equivalent to a second-order cone programming (SOCP) problem [58]. Hence, the computational complexity is given by [58], [59]

$$\mathcal{O}\left(k_{\text{qcqp}}^{0.5}\left(m_{\text{qcqp}}^3 + m_{\text{qcqp}}^2 \sum_{i=1}^{k_{\text{qcqp}}} n_{i,\text{qcqp}} + \sum_{i=1}^{k_{\text{qcqp}}} n_{i,\text{qcqp}}^2\right)\right) \quad (37)$$

where $k_{\text{qcqp}}$ denotes the number of quadratic constraints, $m_{\text{qcqp}}$ denotes the dimension of the optimization variable, and $n_{i,\text{qcqp}}$ denotes the dimension of the $i$-th quadratic constraint. As for the QCQP sub-problem in (22), we have $k_{\text{qcqp}} = N_{mask} + N_{ris}$, where $N_{mask}$ is the number of reradiation mask constraints, $m_{\text{qcqp}} = N_{ris}$, and $n_{i,\text{qcqp}} = N_{ris}, \forall i \in \{1, \ldots, N_{ris} + N_{mask}\}$. Therefore, the complexity of solving the sub-problem in (22) is $\mathcal{O}\left((N_{ris})^3(N_{mask} + N_{ris})^{1.5}\right)$. As for the QCQP sub-problem in (29), we have $k_{\text{qcqp}} = 1$, $m_{\text{qcqp}} = N_t N_{r1}$, and $n_{1,\text{qcqp}} = N_t N_{r1}$. Therefore, the complexity of solving the sub-problem in (29) is $\mathcal{O}\left(N_t^3 N_{r_1}^3\right)$. Since the sub-problems (29) and (34) have the same structure, the computational complexity of (34) is $\mathcal{O}\left(N_t^3 N_{r_2}^3\right)$ as well. Hence, the overall complexity of our proposed algorithm assuming a solver tolerance $\varepsilon = 2.22 \times 10^{-16}$ in CVX is $\mathcal{O}\left(\left(N_{ris}^3(N_{mask} + N_{ris})^{1.5} + N_t^3 N_{r_1}^3 + N_t^3 N_{r_2}^3\right) \log_2(1/\varepsilon)\right) + \mathcal{O}(N_{r_1}^2 N_t + N_{r_1}^3 + N_{r_2}^2 N_t + N_{r_2}^3)$.

## V. Optimization Algorithm for UADP

In this section, we present a greedy search algorithm to solve the optimization problem in (14) when the phase shifts are discrete, i.e., $\theta_n \in \mathcal{F}_2$ [6]. In this case, the problem is challenging due to the non-convex nature of the constraint in (14b). Existing approaches often rely on relaxation methods that first treat the phase shift $\theta_n$ as a continuous variable and then discretize it by projecting onto discrete intervals, such as [23]

$$\theta_n^{\star(\mathcal{F}_2)} = \arg\min_{\theta \in \left\{0, \frac{2\pi}{L}, \cdots, \frac{(L-1)2\pi}{L}\right\}} \left|\theta_n^{\star(\mathcal{F}_1)} - \theta\right| \quad (38)$$

where $\theta_n^{\star(\mathcal{F}_1)}$ represents the optimal reflection coefficient obtained by solving (22), which is a continuous value. However, applying this method alone is insufficient for our problem due to the need to satisfy the reradiation mask constraint in (14d).

The optimal precoders matrices $\mathbf{F}_1$ and $\mathbf{F}_2$, derived by solving (29) and (34), remain unchanged. Therefore, we can still use the same approach to update their values. In what follows, we focus on optimizing $\boldsymbol{\theta}$. Moving from (22), the optimization problem is given by

$$\max_{\boldsymbol{\theta}} \quad \min\left\{-\boldsymbol{\theta}^H \mathbf{E}_1 \boldsymbol{\theta} + 2\operatorname{Re}\left(\boldsymbol{\theta}^H \mathbf{b}_1\right) + c_1,\right.$$
$$\left. -\boldsymbol{\theta}^H \mathbf{E}_2 \boldsymbol{\theta} + 2\operatorname{Re}\left(\boldsymbol{\theta}^H \mathbf{b}_2\right) + c_2\right\} \quad (39a)$$
$$\text{s.t.} \quad \theta_n \in \mathcal{F}_2 \quad \forall n \in \{1, \cdots, N_{ris}\} \quad (39b)$$
$$\Pr(\boldsymbol{\theta}, \mathbf{F}_1, \mathbf{F}_2, \theta^{(ob)}) \leq \rho \quad \theta^{(ob)} \in \mathcal{A} \quad (39c)$$

We propose a greedy search algorithm, described in Algorithm 2, to effectively solve this problem. The algorithm begins by initializing $\boldsymbol{\theta}$ with a continuous value $\boldsymbol{\theta}_0$ that meets the unit amplitude constraint inherent in (39b) and the reradiated mask constraint in (39c). This initialization is performed by normalizing $\boldsymbol{\theta}_0$ to ensure the unit amplitude requirement in (14b), i.e., $\theta_{0,n} = \theta_{0,n}/|\theta_{0,n}|, \forall n$. To ensure that (14d) is fulfilled for $\theta^{(ob)} \in \mathcal{A}$, we introduce the maximum reradiated power, based on (20), as follows:

$$\Pr^{(max)} = \max_{\theta^{(ob)} \in \mathcal{A}} \left(\boldsymbol{\theta}_0^H \left(\mathbf{Q}^{(ob)} \odot \mathbf{T}_1 + \mathbf{Q}^{(ob)} \odot \mathbf{T}_2\right) \boldsymbol{\theta}_0\right) \quad (40)$$

[6] We use a greedy search since it has low complexity and is more suitable for the discrete phase shifts in practical RIS hardware, where the solution space is already limited by the phase quantization.

**Algorithm 2** Greedy search algorithm for solving problem (39)

1: **Input:** $t = 0$, $\epsilon > 0$, $\boldsymbol{\theta}^{(0)} = \boldsymbol{\theta}_0$ to fulfill unit amplitude and (39c);
2: Compute $f^{(0)} = f(\boldsymbol{\theta}^{(0)})$ from (39a);
3: Set an indicator variable $\mathrm{T}^{ind} = 0$;
4: **repeat**
5:   **for** $n = 1, \ldots, N_{ris}$ **do**
6:     Initialize a max value $f^{max} = \inf$;
7:     **for** $l = 1, \ldots, L$ **do**
8:       Compute $\theta_n^{mid} = \exp(j\frac{2\pi}{L}(l-1))$;
9:       $\boldsymbol{\theta}^{mid} \leftarrow [\theta_1, \ldots, \theta_n^{mid}, \ldots, \theta_{N_{ris}}]^T$;
10:       Compute the indicator function $\mathrm{T}(\boldsymbol{\theta}^{mid})$;
11:       **if** $\mathrm{T}(\boldsymbol{\theta}^{mid}) = 1$ **then**
12:         **if** $f(\boldsymbol{\theta}^{mid}) > f^{max}$ **then**
13:           Update $\theta_n \leftarrow \theta_n^{mid}$ and $f^{max} \leftarrow f(\boldsymbol{\theta}^{mid})$;
14:         **end if**
15:       **end if**
16:     **end for**
17:   **end for**
18:   $f^{(t+1)} \leftarrow f(\boldsymbol{\theta})$, $\mathrm{T}^{ind} \leftarrow \mathrm{T}(\boldsymbol{\theta})$
19:   $t \leftarrow t + 1$
20: **until** $|f^{(t)} - f^{(t-1)}| < \epsilon$ and $\mathrm{T}^{ind} = 1$

where we note that $\mathbf{Q}^{(ob)}$ is dependent on $\theta^{(ob)} \in \mathcal{A}$. If $\Pr^{(max)} \leq \rho$, the initial value remains the same where $\tilde{\boldsymbol{\theta}}_0 = P_\theta(\boldsymbol{\theta}_0) = \boldsymbol{\theta}_0$. Otherwise, we compute the normalized reflection coefficient $\tilde{\boldsymbol{\theta}}_0 = a\boldsymbol{\theta}_0$, where $a$ is obtained such that

$$\max_{\theta^{(ob)} \in \mathcal{A}} \left(\tilde{\boldsymbol{\theta}}_0^H \left(\mathbf{Q}^{(ob)} \odot \mathbf{T}_1 + \mathbf{Q}^{(ob)} \odot \mathbf{T}_2\right) \tilde{\boldsymbol{\theta}}_0\right) \quad (41)$$
$$= a^2 \max_{\theta^{(ob)} \in \mathcal{A}} \left(\boldsymbol{\theta}_0^H \left(\mathbf{Q}^{(ob)} \odot \mathbf{T}_1 + \mathbf{Q}^{(ob)} \odot \mathbf{T}_2\right) \boldsymbol{\theta}_0\right)$$
$$= \rho$$

This results in the following initial value, and in the following projection function $P_\theta(\boldsymbol{\theta}_0)$

$$\boldsymbol{\theta}_0 = P_\theta(\boldsymbol{\theta}_0)$$
$$= \boldsymbol{\theta}_0 \sqrt{\rho / \max_{\theta^{(ob)} \in \mathcal{A}} \left(\boldsymbol{\theta}_0^H \left(\mathbf{Q}^{(ob)} \odot \mathbf{T}_1 + \mathbf{Q}^{(ob)} \odot \mathbf{T}_2\right) \boldsymbol{\theta}_0\right)}$$
(42)

The algorithm then iteratively optimizes each element of $\boldsymbol{\theta}$ by selecting values from the discrete set

$$\mathcal{F}_2 = \left\{\theta_n \Big| \theta_n = e^{j\frac{2\pi}{L}(l-1)}, \ l \in \{1, \ldots, L\}\right\} \quad (43)$$

Next, a greedy approach is applied. For each phase shift element, i.e., $\theta_n$, $n = 1, 2, \ldots, N_{ris}$, the algorithm evaluates all possible values in (43), selecting the one that maximizes the objective function in (39a), denoted as $f$, while satisfying the mask constraint. An indicator function $\mathrm{T}(\boldsymbol{\theta})$ is employed to verify whether the current $\boldsymbol{\theta}$ violates the mask constraint. The indicator function is defined as

$$\mathrm{T}(\boldsymbol{\theta}) = \begin{cases} 1, & \text{if } \Pr(\boldsymbol{\theta}, \mathbf{F}_1, \mathbf{F}_2, \theta^{(ob)}) \leq \rho, \ \theta^{(ob)} \in \mathcal{A} \\ 0, & \text{otherwise} \end{cases}$$

The process is repeated until convergence, defined by the change in the objective function becoming smaller than a predefined threshold $\epsilon$, and the indicator function confirming that the mask constraint is satisfied. This greedy approach ensures that the algorithm converges to a solution that maximizes the minimum of the achievable rates while fulfilling all the constraints. In this regard, it is important to note that (41) is applied only during the initialization phase. This is to accelerate the search process by generating initial candidate values that are more likely to satisfy the mask constraint.

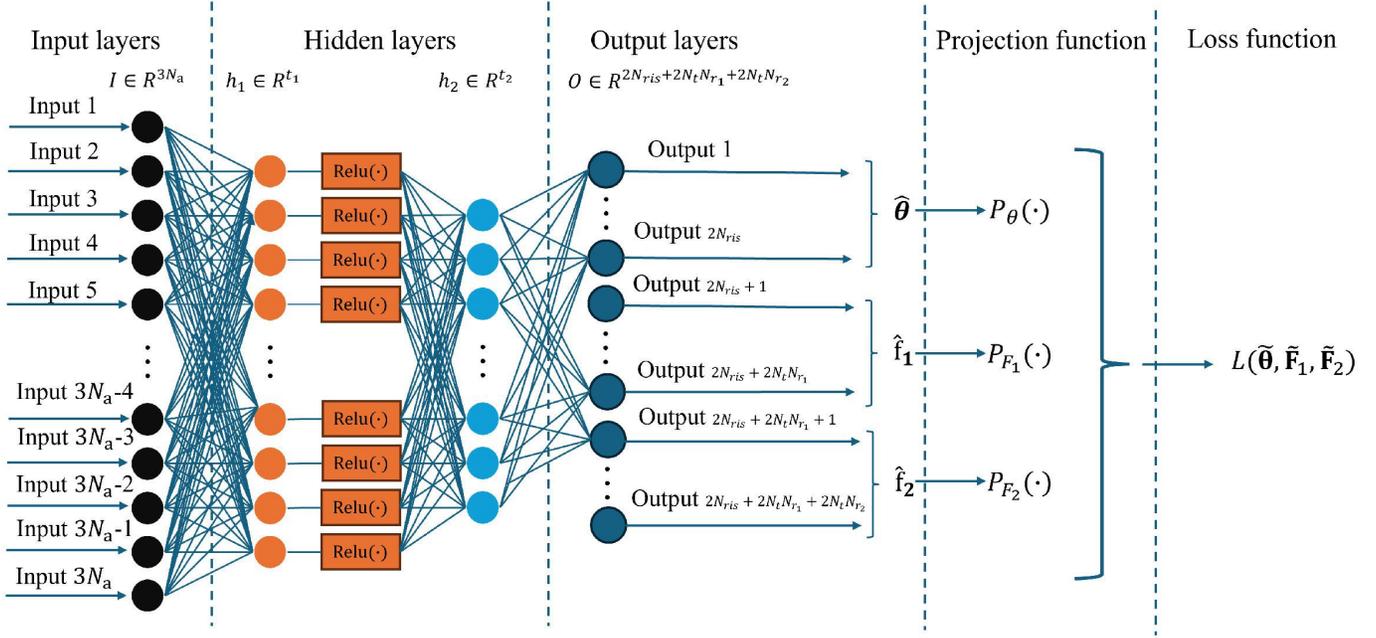

Fig. 2. Considered neural network architecture

**Complexity** – The complexity of the algorithm can be expressed as $\mathcal{O}(M_{\text{iter}} N_{\text{ris}} L)$ [7], where $M_{\text{iter}}$ is the number of iterations required for convergence, $N_{\text{ris}}$ is the number of RIS elements, and $L$ is the number of phase shift levels, i.e., the number of possible phase shift values.

**Convergence** – The algorithm evaluates all possible values of each phase shift element $\theta_n$ from the discrete set $\mathcal{F}_2$, selecting the one that maximizes the objective function $f$, while satisfying the mask constraint. This guarantees that the objective function either increases or remains constant, ensuring a monotonic improvement at each step. Since the search space is finite, the algorithm will eventually try all possible updates and stop when no further improvement is possible, which is determined by the threshold $\epsilon$. Thus, the algorithm is guaranteed to converge.

## VI. MODEL-DRIVEN NEURAL NETWORK OPTIMIZATION FOR UACP

In this section, we present an approach for RIS optimization that improves the performance of the algorithm introduced in the previous section. Specifically, we provide a more efficient implementation by incorporating a neural network. Unlike conventional data-driven neural networks, our approach is model-driven. Specifically, our approach leverages the mathematical structure of the original optimization problem, introduced in Section IV, to guide the network design and the training process. The network is trained using a tailored loss function that directly reflects the objective function of the optimization problem, which reduces the implementation complexity. The details of this approach are presented next. The proposed neural network architecture takes as inputs the angle of incidence $\theta_{inc}$ and the desired angles of reflection $\theta_{ref1}$ and $\theta_{ref2}$, and returns the concatenated vector of reflection coefficients of the RIS, along with the vectorized precoders $\mathbf{f}_1$ and $\mathbf{f}_2$ in (25) and (32), respectively.

The proposed neural network architecture is shown in Fig. 2. Given that the reflection coefficients $\boldsymbol{\theta}$ and the vectorized precoders $\mathbf{f}_1$ and $\mathbf{f}_2$ are complex-valued vectors, the outputs are presented in terms of real and imaginary parts of $\boldsymbol{\theta}$, $\mathbf{f}_1$ and $\mathbf{f}_2$, as follows:

$$\hat{\boldsymbol{\theta}} = \begin{bmatrix} \mathcal{R}(\boldsymbol{\theta}) \\ \mathcal{J}(\boldsymbol{\theta}) \end{bmatrix} \in \mathbb{R}^{2N_{\text{ris}} \times 1} \tag{44}$$

$$\hat{\mathbf{f}}_1 = \begin{bmatrix} \mathcal{R}(\mathbf{f}_1) \\ \mathcal{J}(\mathbf{f}_1) \end{bmatrix} \in \mathbb{R}^{2N_t N_{r_1} \times 1} \tag{45}$$

$$\hat{\mathbf{f}}_2 = \begin{bmatrix} \mathcal{R}(\mathbf{f}_2) \\ \mathcal{J}(\mathbf{f}_2) \end{bmatrix} \in \mathbb{R}^{2N_t N_{r_2} \times 1} \tag{46}$$

These three vectors are then concatenated as an output of the neural network where

$$\hat{\mathbf{y}} = \begin{bmatrix} \hat{\boldsymbol{\theta}} \\ \hat{\mathbf{f}}_1 \\ \hat{\mathbf{f}}_2 \end{bmatrix} \in \mathbb{R}^{(2N_{\text{ris}} + 2N_t N_{r_1} + 2N_t N_{r_2}) \times 1} \tag{47}$$

To ensure that the neural network implementation provides accurate results, we propose employing a modified one-hot encoding method [60]. The one-hot encoding method is motivated by the challenge that neural networks may encounter in differentiating between angles that vary by only a small fraction of a degree. Directly feeding raw continuous angle values makes it difficult for the network to distinguish small angular differences, which can significantly impact the beamforming accuracy, as shown later in the numerical results. While traditional one-hot encoding works well for discrete values, it is not suitable for representing continuous angles. The proposed one-hot encoding method enlarges the differences between similar angles. Unlike the conventional one-hot encoding, which represents inputs strictly as zeros and ones, our method incorporates the fractional components of the angles.

---
[7]Note that the proposed greedy algorithm is different from the exhaustive search, which evaluates all possible combinations of RIS phase shifts, resulting in an exponential complexity of $\mathcal{O}(L^{N_{ris}})$. Hence, the greedy algorithm has a significant lower complexity compared to the exhaustive search.

Specifically, the one-hot function applied to the angle of incidence is defined as follows

$$\boldsymbol{y}_{inc} = \text{one\_hot}(\theta_{inc}) \tag{48}$$

where $\theta_{inc}$ is assumed to lie within the range $\theta_{low} \leq \theta_{inc} \leq \theta_{high}$. For simplicity, we assume that $\theta_{low}$ and $\theta_{high}$ are integer values. The desired angular resolution is denoted by $\mu$. Therefore, there exist $N_a = (\theta_{high} - \theta_{low})/\mu + 1$ possible values of interest for $\theta_{inc}$. For convenience, we introduce the function $\text{frac}(\cdot)$, which returns the fractional part of a real number, and the function $\text{int}(\cdot)$, which returns the integer part of a real number. Then, $\boldsymbol{y}_{inc} \in \mathbb{R}^{N_a \times 1}$ is expressed as

$$\boldsymbol{y}_{inc}^T = [\ \underbrace{\boldsymbol{0}_1^T}_{\mathbb{R}^{1 \times \text{int}((\theta_{inc}-\theta_{low})/\mu)}}, 1 + \text{frac}(\theta_{inc}) - \mu, \boldsymbol{0}_1^T] \tag{49}$$

if $\text{frac}(\theta_{inc}) > \mu$ and

$$\boldsymbol{y}_{inc}^T = [\ \underbrace{\boldsymbol{0}_1^T}_{\mathbb{R}^{1 \times \text{int}((\theta_{inc}-\theta_{low})/\mu)}}, \text{frac}(\theta_{inc}), \boldsymbol{0}_1^T] \tag{50}$$

otherwise, where $\boldsymbol{0}_1^T$ is an 1-dimensional zero vector. The same encoding method can be applied to the angles $\theta_{ref1}$ and $\theta_{ref2}$ to obtain $\boldsymbol{y}_{ref1} \in \mathbb{R}^{N_a \times 1}$ and $\boldsymbol{y}_{ref2} \in \mathbb{R}^{N_a \times 1}$, respectively. The input to the neural network in Fig. 2 is then formed by stacking $\boldsymbol{y}_{inc}$, $\boldsymbol{y}_{ref1}$ and $\boldsymbol{y}_{ref}$, as follows:

$$\boldsymbol{x} = [\boldsymbol{y}_{inc}^T, \boldsymbol{y}_{ref1}^T, \boldsymbol{y}_{ref2}^T]^T \in \mathbb{R}^{3N_a \times 1} \tag{51}$$

The considered neural network architecture consists of two hidden layers. The numbers of neurons in the first and second hidden layers are denoted by $t_1$ and $t_2$, respectively. These values are determined through a trial-and-error approach. The vectors $\boldsymbol{s}_1 \in \mathbb{R}^{t_1 \times 1}$ and $\boldsymbol{s}_2 \in \mathbb{R}^{t_2 \times 1}$ denote the output data from the first and second hidden layers, respectively. The weight matrix and bias vector for the first hidden layer are denoted by $\mathbf{W}^{(1)} \in \mathbb{R}^{t_1 \times 3N_a}$ and $\mathbf{b}^{(1)} \in \mathbb{R}^{t_1 \times 1}$; for the second hidden layer, the weight matrix and bias vector are given by $\mathbf{W}^{(2)} \in \mathbb{R}^{t_2 \times t_1}$ and $\mathbf{b}^{(2)} \in \mathbb{R}^{t_2 \times 1}$; in the output layer, the weight matrix and bias vector are $\mathbf{W}^{(3)} \in \mathbb{R}^{(2N_{ris}+2N_tN_{r_1}+2N_tN_{r_2}) \times t_2}$ and $\mathbf{b}^{(3)} \in \mathbb{R}^{(2N_{ris}+2N_tN_{r_1}+2N_tN_{r_2}) \times 1}$. Also, a Relu activation function is applied to the hidden layers. In mathematical terms, this is expressed as

$$\boldsymbol{s}_1 = \text{Relu}(\mathbf{W}^{(1)}\boldsymbol{x} + \mathbf{b}^{(1)}) \tag{52}$$
$$\boldsymbol{s}_2 = \mathbf{W}^{(2)}\boldsymbol{s}_1 + \mathbf{b}^{(2)} \tag{53}$$
$$\hat{\boldsymbol{y}} = \mathbf{W}^{(3)}\boldsymbol{s}_2 + \mathbf{b}^{(3)} \tag{54}$$

We first extract $\hat{\boldsymbol{\theta}}$, $\hat{\boldsymbol{f}_1}$ and $\hat{\boldsymbol{f}_2}$ from the output $\hat{\boldsymbol{y}}$, where

$$\hat{\boldsymbol{\theta}} = \hat{\boldsymbol{y}}[1:2N_{ris}] \tag{55}$$
$$\hat{\mathbf{f}}_1 = \hat{\boldsymbol{y}}[2N_{ris}+1:2N_{ris}+2N_tN_{r_1}] \tag{56}$$
$$\hat{\mathbf{f}}_2 = \hat{\boldsymbol{y}}[2N_{ris}+2N_tN_{r_1}+1:2(N_{ris}+N_tN_{r_1}+N_tN_{r_2})] \tag{57}$$

Since $\hat{\boldsymbol{\theta}}$, $\hat{\mathbf{f}}_1$ and $\hat{\mathbf{f}}_2$ are real values, we need to retrieve the complex reflection coefficient $\bar{\boldsymbol{\theta}}$ and the complex precoder vectors $\bar{\mathbf{f}}_1$ and $\bar{\mathbf{f}}_2$ where

$$\mathcal{R}(\bar{\boldsymbol{\theta}}) = \hat{\boldsymbol{\theta}}[1:N_{ris}+1,1] \in \mathbb{R}^{N_{ris} \times 1} \tag{58}$$
$$\mathcal{J}(\bar{\boldsymbol{\theta}}) = \hat{\boldsymbol{\theta}}[N_{ris}+1:2N_{ris},1] \in \mathbb{R}^{N_{ris} \times 1} \tag{59}$$

$$\mathcal{R}(\bar{\mathbf{f}}_1) = \hat{\mathbf{f}}_1[1:N_tN_{r_1},1] \in \mathbb{R}^{N_tN_{r_1} \times 1} \tag{60}$$
$$\mathcal{J}(\bar{\mathbf{f}}_1) = \hat{\mathbf{f}}_1[N_tN_{r_1}+1:2N_tN_{r_1},1] \in \mathbb{R}^{N_tN_{r_1} \times 1} \tag{61}$$

$$\mathcal{R}(\bar{\mathbf{f}}_2) = \hat{\mathbf{f}}_2[1:N_tN_{r_2},1] \in \mathbb{R}^{N_tN_{r_2} \times 1} \tag{62}$$
$$\mathcal{J}(\bar{\mathbf{f}}_2) = \hat{\mathbf{f}}_2[N_tN_{r_2}+1:2N_tN_{r_2},1] \in \mathbb{R}^{N_tN_{r_2} \times 1} \tag{63}$$

The precoding vectors $\bar{\mathbf{f}}_1 \in \mathbb{C}^{N_tN_{r_1} \times 1}$ and $\bar{\mathbf{f}}_2 \in \mathbb{C}^{N_tN_{r_2} \times 1}$ are then reshaped into their corresponding matrix forms, $\bar{\mathbf{F}}_1 \in \mathbb{C}^{N_t \times N_{r_1}}$ and $\bar{\mathbf{F}}_2 \in \mathbb{C}^{N_t \times N_{r_2}}$, following the method in (30) and (35), respectively, where

$$\bar{\mathbf{F}}_1 = [\bar{\mathbf{f}}_1(1:N_t), \cdots, \bar{\mathbf{f}}_1((N_{r_1}-1)N_t+1:N_tN_{r_1})] \tag{64}$$
$$\bar{\mathbf{F}}_2 = [\bar{\mathbf{f}}_2(1:N_t), \cdots, \bar{\mathbf{f}}_2((N_{r_2}-1)N_t+1:N_tN_{r_2})] \tag{65}$$

It is necessary that the precoding matrices $\bar{\mathbf{F}}_1$ and $\bar{\mathbf{F}}_2$ generated by the neural network satisfy the transmit power constraint in (14c). To ensure this, the precoding matrices need to be projected to meet the constraint. The total transmit power is computed as

$$P_t = \sum_{i=1}^{2} \text{Tr}(\bar{\mathbf{F}}_i \bar{\mathbf{F}}_i^H) \tag{66}$$

If the total transmit power satisfies the constraint, no projection is necessary. However, if the constraint is violated, the precoding matrices are scaled as follows

$$\tilde{\mathbf{F}}_i = \alpha \bar{\mathbf{F}}_i, \quad i \in \{1, 2\} \tag{67}$$

where $\alpha$ is the scaling factor defined by

$$\alpha = \sqrt{P_{max}/P_t} \tag{68}$$

Additionally, the solution provided by the neural network must satisfy the unit amplitude constraint in (14b) and the reradiation mask constraint in (14d) over the specified set of angles $\mathcal{A}$. To achieve this, we first normalize $\bar{\boldsymbol{\theta}}$ by applying $\bar{\theta}_n = \bar{\theta}_n/|\bar{\theta}_n|, \forall n$ to meet the unit amplitude requirement in (14b). This normalization is always applied to ensure the unit amplitude constraint is met, regardless of whether the mask constraint is satisfied. However, when the normalized vector $\bar{\boldsymbol{\theta}}$ violates the reradiation mask constraint, a scaling factor is applied using the projection function $\tilde{\boldsymbol{\theta}} = P_\theta(\bar{\boldsymbol{\theta}})$ as described in (42). This ensures the mask constraint is satisfied, even if the unit amplitude condition is relaxed. This tradeoff is necessary in those cases where both constraints cannot be satisfied simultaneously, which depends on the maximum value of the reradiation mask.

The proposed model-based neural network is based on the idea that training data and labels can be avoided by appropriately choosing the loss function derived from a reliable model, while leveraging the computing capabilities of the neural network architecture and its network structure. Based on Section IV, we consider the following loss function for optimizing the biases and the weights of the proposed network architecture:

$$L(\tilde{\boldsymbol{\theta}}, \tilde{\mathbf{F}}_1, \tilde{\mathbf{F}}_2) = \min_{\tilde{\boldsymbol{\theta}}, \tilde{\mathbf{F}}_1, \tilde{\mathbf{F}}_2} -\min\left\{R_1(\tilde{\boldsymbol{\theta}}, \tilde{\mathbf{F}}_1, \tilde{\mathbf{F}}_2), R_2(\tilde{\boldsymbol{\theta}}, \tilde{\mathbf{F}}_1, \tilde{\mathbf{F}}_2)\right\} \quad (69)$$

Thanks to (52)-(57), the loss function in (69) depends explicitly on the weights and biases, which are optimized using the Adagrad adaptive gradient algorithm. In this approach, the neural network architecture is predetermined (number of layers and number of neurons in each layer), and the weights and biases are obtained from the model specified by the loss function in (69).

## VII. SIMULATION RESULTS

In this section, we validate the effectiveness of the proposed algorithms. As depicted in Fig. 1, we consider a three-dimensional coordinate system. In this setup, the midpoint of the RIS is positioned at the origin $(0, 0, 0)$, the midpoint of transmitter is located at $(-d_{ris}, 0, l_t)$, and the midpoints of the receivers $r_1$ and $r_2$ are located at $(D - d_{ris}, 0, l_{r_1})$ and $(D - d_{ris}, 0, l_{r_2})$, respectively. Based on the network geometry, we have $l_t = d_{ris}/\tan(\theta_{inc})$, $l_{r_1} = (D - d_{ris})/\cos(\theta_{ref1})$ and $l_{r_2} = (D - d_{ris})/\cos(\theta_{ref2})$, respectively. The angles are set such that $\theta_{low} = 10°$ and $\theta_{high} = 60°$. Thus, the angle of incidence $\theta_{inc}$, and the angles of reflection $\theta_{ref1}$ and $\theta_{ref2}$ lie within the range $[10, 60]$ degrees. We assume an angular resolution of $\mu = 0.5°$, so the input layer of the neural network has size $d = 3$ and $N_a = 303$. We set the convergence tolerance to $\epsilon = 10^{-8}$ and $L = 4$ are the phase shift levels. The simulation parameters are provided in Table I, unless stated otherwise, and the neural network configuration is given in Table II. Specifically, the simulation parameters are chosen based on trial runs to ensure stability and fast convergence. Also, the reradiation mask, represented by the set $\mathcal{A}$, is chosen as follows:

$$\mathcal{A} \in [-89, \theta_{ref1} - 15] \cup [\theta_{ref2} + 15, 89] \quad (70)$$

if $\theta_{ref2} - \theta_{ref1} \leq 20$ and

$$\mathcal{A} \in [-89, \theta_{ref1} - 15] \cup [\theta_{ref2} + 15, 89] \\ \cup [\theta_{ref1} + 10, \theta_{ref2} - 10] \quad (71)$$

if $\theta_{ref2} - \theta_{ref1} > 20$.

For comparison, we consider four different schemes: (1) AO for UACP without imposing the mask constraint in (14d); (2) AO for UACP imposing the mask constraint in (14d); (3) model-based neural network for UACP imposing the mask constraint in (14d); and (4) greedy search for UADP imposing the mask constraint in (14d).

Table III presents the running time of the four considered schemes, which is obtained using a PC equipped with an

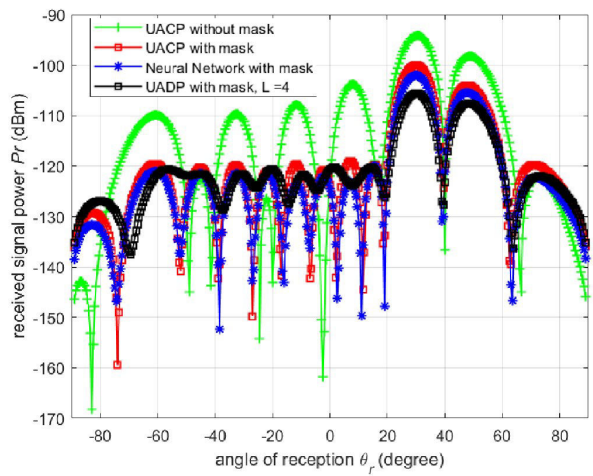

(a) Received power as a function of the angle of observation

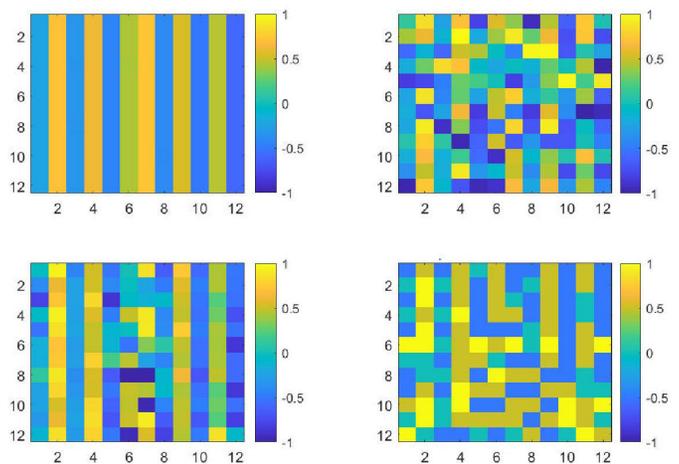

(b) Color map of the RIS phase shifts: **UACP** without mask (top left), **UACP** with mask (top right), neural network with mask (bottom left), **UADP** with mask (bottom right)

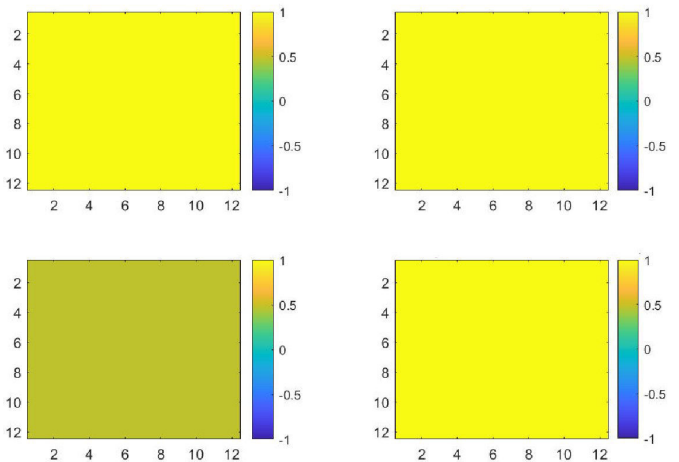

(c) Color map of the RIS amplitudes: **UACP** without mask (top left), **UACP** with mask (top right), neural network with mask (bottom left), **UADP** with mask (bottom right)

Fig. 3. Results for the setup $\theta_{inc} = 20°$, $\theta_{ref1} = 30°$, $\theta_{ref2} = 50°$

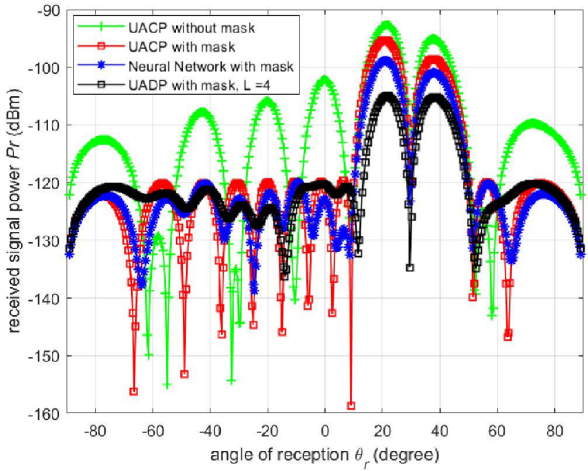

(a) Received power as a function of the angle of observation

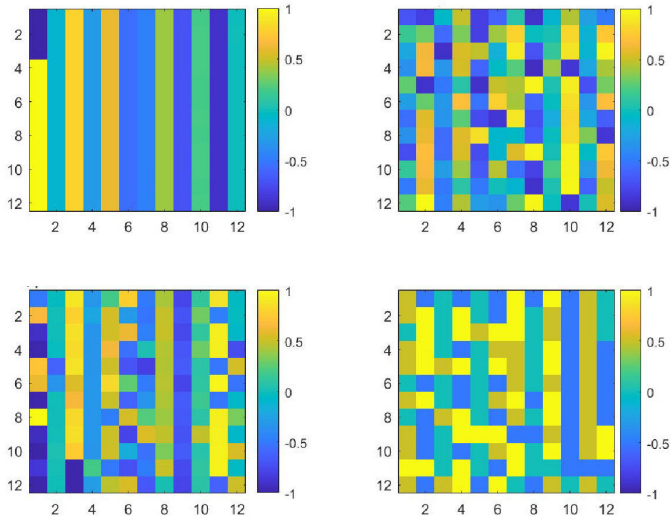

(b) Color map of the RIS phase shifts: **UACP** without mask (top left), **UACP** with mask (top right), neural network with mask (bottom left), **UADP** with mask (bottom right)

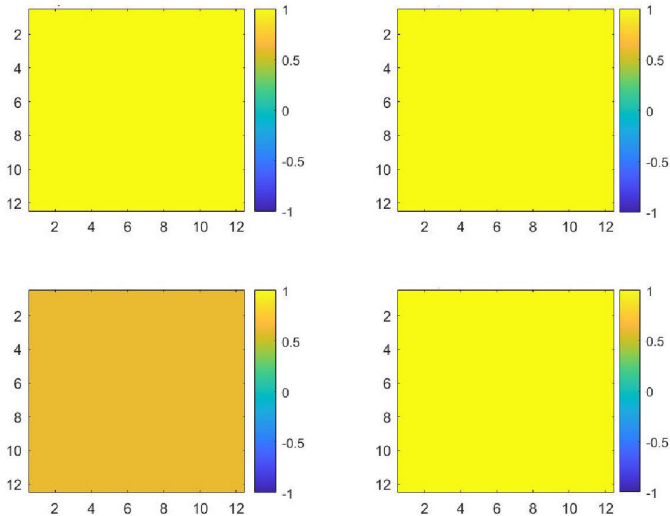

(c) Color map of the RIS amplitudes: **UACP** without mask (top left), **UACP** with mask (top right), neural network with mask (bottom left), **UADP** with mask (bottom right)

Fig. 4. Results for the setup $\theta_{inc} = 25°$, $\theta_{ref1} = 20°$, $\theta_{ref2} = 40°$

TABLE I
SIMULATION PARAMETERS

| Parameters | Values |
|---|---|
| Frequency | 27 GHz |
| $K_r$ | $10^5$ |
| $D$ | 100 m |
| $d_{ris}$ | 20 m |
| $\rho$ | -120 dBm |
| $P_{max}$ | 2 Watt |
| $\sigma^2$ | -120 dB Watt |
| $N_t$ | 8 |
| $N_{r_1}, N_{r_2}$ | 2 |
| $N_{ris}$ | $12 \times 12$ |

TABLE II
NEURAL NETWORK PARAMETERS

| Parameter | Setup |
|---|---|
| optimizer | adaptive gradient |
| learning rate | 0.0005 |
| iterations | 100 |
| $t_1$ | 1024 |
| $t_2$ | 512 |
| $\mu$ | 0.5 |

Intel(R) Core(TM) i9-10900 CPU and 32.0 GB RAM [8]. The results indicate that imposing the reradiation mask constraints significantly increases the running time. Importantly, the proposed model-based neural network achieves a reduced running time compared to the AO method for the UACP case. In contrast, due to the application of the greedy search algorithm for the UADP case, this latter scheme has the highest running time.

TABLE III
RUNNING TIME COMPARISON

| Schemes | Running time (seconds) |
|---|---|
| UACP without mask | 40.69 |
| UACP with mask | 731.4 |
| Neural Network with mask | 120.98 |
| UADP with mask | 859.27 |

In Figs. 3 and 4, we analyze the power reradiated by the RIS for any angle of observation, considering two different case studies for the angles of incidence and desired angles of reflection. Also, we present the corresponding amplitude and phase of the reflection coefficients of the RIS elements. When the reradiation mask constraint is not imposed, we observe that the side lobes can be quite high. However, when the reradiation mask is imposed, the constraint is satisfied, and the main lobes remain directed toward the desired directions. By imposing the reradiation mask, the gain of the main lobes is typically reduced, but not significantly. The neural network implementation produces good results comparable to the AO optimization framework for the UACP case, with only a minor decrease in the gain of the main lobes. By analyzing the amplitude of the reflection coefficient, we observe that $|\tilde{\boldsymbol{\theta}}_i| \approx 1$ in the case studies considered. However, when the reradiation mask is imposed, some RIS elements have $|\tilde{\boldsymbol{\theta}}_i| < 1$. This reduction is due to the necessity of fulfilling the reradiation mask constraint, which imposes some fundamental limits on the reradiation efficiency of the RIS. In the UADP case, the main lobes are still directed toward the desired directions, with a small reduction in the gain of the main lobes due to the phase shifts being restricted to $L = 4$ possible discrete values.

---

[8]Note that the proposed algorithms are intended to run at the BS, where hardware accelerators like graphics processing units (GPUs) or field-programmable gate arrays (FPGAs) can significantly reduce the computation time.

<pre>

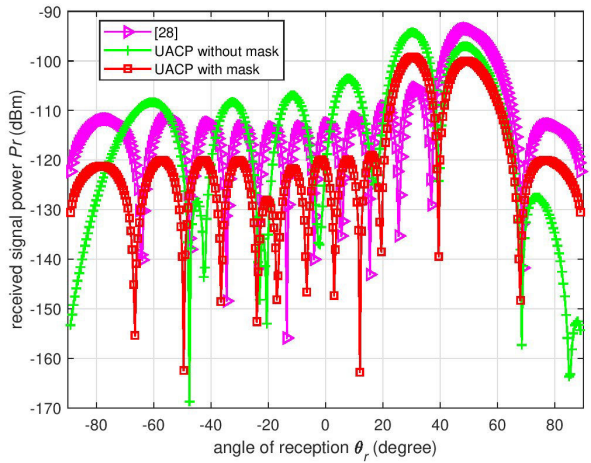

Fig. 5. Received power as a function of the angle of observation for the setup $\theta_{inc} = 20°$, $\theta_{ref2} = 30°$, $\theta_{ref2} = 50°$, $N_t = 8$ and $N_{r_1} = N_{r_2} = 2$

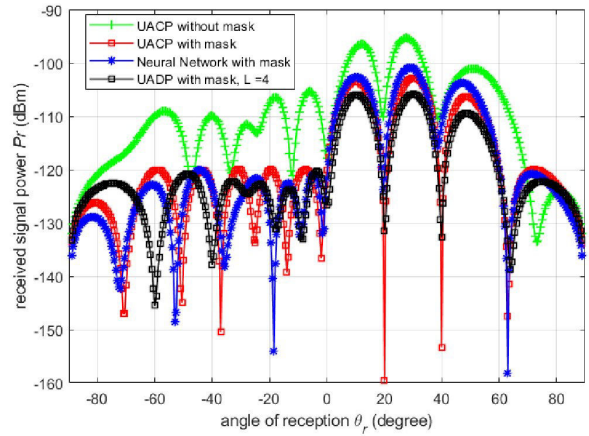

Fig. 7. Received power as a function of the angle of observation for the setup $\theta_{inc} = 20°$, $\theta_{ref1} = 10°$, $\theta_{ref2} = 30°$, $\theta_{ref3} = 50°$, $N_t = 8$ and $N_{r_1} = N_{r_2} = N_{r_3} = 2$

Fig 5 compares the beam patterns of our proposed UACP methods, with and without mask, to the algorithm in [28]. We can see that the method in [28] does not ensure proper alignment of the main lobes with both target directions, leading to suboptimal beamforming and significant side lobes in undesired directions. In contrast, the UACP method with the mask achieves a better signal power distribution across the two receivers and reduced interference towards unwanted directions.

receivers located in distinct directions as shown in Fig. 7. The results show that the proposed methods are capable of handling multi-beam beamforming where the received signal power is well distributed across the three beams, which ensures balanced performance across all targeted directions, while the side lobes are effectively suppressed in all other undesired directions.

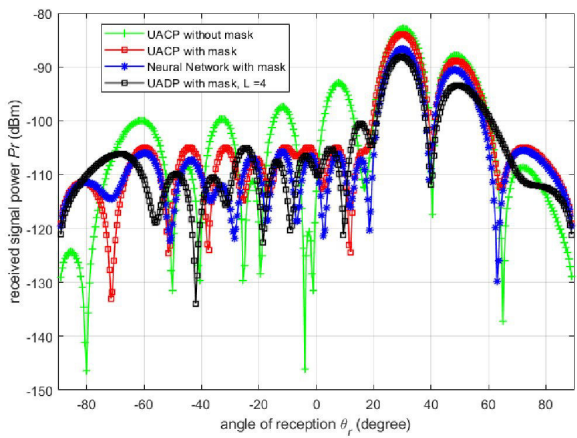

Fig. 6. Received power as a function of the angle of observation for the setup $\theta_{inc} = 20°$, $\theta_{ref1} = 30°$, $\theta_{ref2} = 50°$, $N_t = 64$ and $N_{r_1} = N_{r_2} = 4$

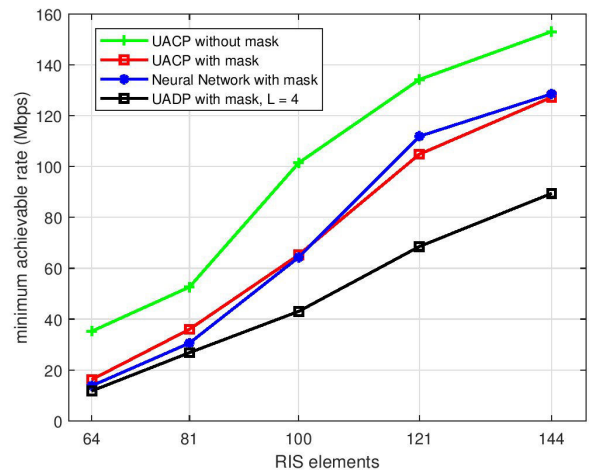

Fig. 8. The minimum achievable rate versus the number of RIS elements for the setup $\theta_{inc} = 20°$, $\theta_{ref1} = 30°$, $\theta_{ref2} = 50°$, $N_t = 8$ and $N_{r_1} = N_{r_2} = 2$

To further evaluate the effectiveness of our proposed methods, we present the results with a larger number of antennas, where $N_t = 64$ and $N_{r_1} = N_{r_2} = 4$ as shown in Fig. 6, where we set $\rho = -105$ dBm. We observe that the main lobes remain directed towards the desired directions while satisfying the reradiation mask constraints. As expected, we also observe an increase in the received signal power, due to the higher array gain provided by the larger number of antennas.

To demonstrate the effectiveness of each proposed method in a multi-beam beamforming scenario, we consider the case where the RIS is configured to direct signals towards three

We illustrate in Fig. 8 the impact of the number of RIS elements on the minimum achievable rate measured in Mbps for the setup $\theta_{inc} = 20°$, $\theta_{ref1} = 30°$, $\theta_{ref2} = 50°$, with $N_t = 8$, $N_{r_1} = N_{r_2} = 2$ and $D = 50$ m. First, we see that the minimum achievable rate increases with the number of RIS elements, which is expected since deploying more RIS elements significantly enhances the communication performance. Among the proposed methods, UACP without mask achieves the highest minimum rate. The neural network-based method and UACP with mask achieves comparable
</pre>

performance, and both show a lower minimum achievable rate compared to UACP without mask. This is expected since the reradiation mask constraint reduces the feasible set of solutions. Additionally, we observe that UADP achieves performance close to that of UACP and the neural network-based methods, with continuous phase shifts, when the number of RIS elements is small. However, as the RIS size increases, the performance gap widens due to the limited phase shift resolution, which increasingly impacts the achievable rate.

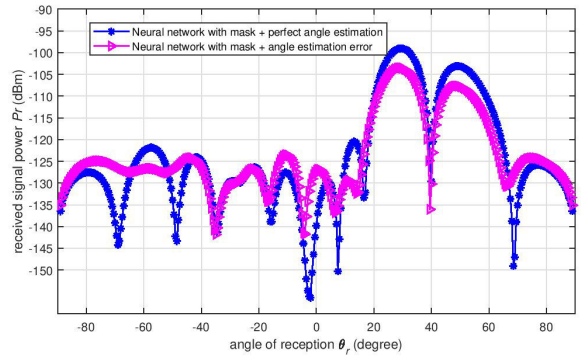

Fig. 10. Beam pattern comparison obtained with the neural network assuming perfect angle estimation and in the presence of estimation errors, for the setup $\theta_{\text{inc}} = 20°$, $\theta_{\text{ref1}} = 30°$, $\theta_{\text{ref2}} = 50°$, with $N_t = 8$ and $N_{r_1} = N_{r_2} = 2$. A small uniform noise in the range $[0, 0.5]$ degrees is added to the input angles to model the estimation error.

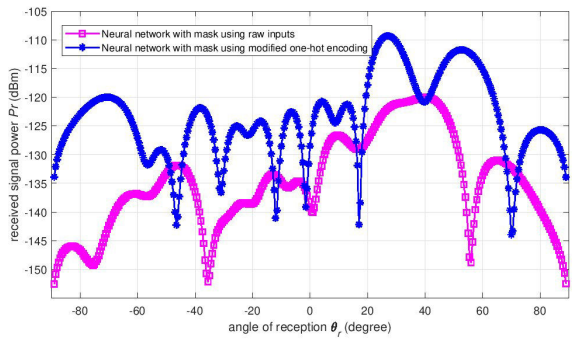

Fig. 9. Beam pattern comparison between the neural network with raw input angles and our proposed modified one-hot encoding for the setup $\theta_{inc} = 20°$, $\theta_{ref1} = 30°$, $\theta_{ref2} = 50°$, $N_t = 8$ and $N_{r_1} = N_{r_2} = 2$

Fig. 9 compares the beamforming performance of the proposed modified one-hot encoding method against the case study where the angles of incidence and reflection are input to the neural network without any encoding (denoted as raw inputs). The results show that our encoding approach enables the neural network to form accurate beams towards the desired directions while effectively suppressing the side lobes in unwanted directions. In contrast, feeding the neural network directly with the angles without any encoding leads to degraded performance, as the neural network cannot capture small differences between closely spaced angles, particularly between the two directions of interests, which results in inaccurate beam patterns. This comparison highlights the effectiveness of the proposed encoding method in improving the beamforming accuracy.

Lastly, to evaluate the robustness of the proposed neural network-based method to the presence of angle estimation errors, we consider a scenario where the input angles are affected by random noise that is uniformly distributed in the range $[0, 0.5]$ degrees. As shown in Fig. 10, despite the presence of estimation errors, the neural network still directs the beams effectively towards the desired directions, with only a slight degradation in the received signal power compared to the case of perfect angle estimation. This shows that our proposed method is robust to small inaccuracies in angle estimation.

From the numerical results obtained, it can be concluded that the proposed methods effectively shape the reradiation pattern of an RIS towards the desired directions while adhering to the specified reradiation mask constraints.

## VIII. CONCLUSION

In this paper, we have addressed the joint design of the transmit precoding matrices and RIS phase shift vector in a two-user RIS-aided MIMO system, focusing on maximizing the minimum achievable rate under transmit power and reradiation mask constraints. We simplified the problem using the Arimoto-Blahut algorithm and solved it through QPQC subproblems using an alternating optimizaion approach.

To improve efficiency, we developed a model-based neural network based on the one-hot encoding for the angles of incidence and reflection, which significantly reduces the execution time while effectively shaping the radiation pattern to meet the specified reradiation mask constraints. We also addressed practical RIS limitations by implementing a greedy search algorithm, assuming discrete-valued phase shifts. Simulation results demonstrated that the proposed methods achieve the desired beam shaping while adhering to the specified reradiation constraints. The neural network approach was shown to offer an efficient alternative to traditional optimization methods. Also, we demonstrated that the use of four discrete values for the phase shifts of the RIS provides a small reduction of the gain, contributing to the practical deployment of RISs in future wireless systems.

## APPENDIX

Let us assume that we have $N$ receivers, each requiring a beam directed towards it, resulting in a total of $N$ beams. The received signal for the receiver $r_i$, $i = \{1, \ldots, N\}$ equipped with $N_{r_i}$ antennas is given by

$$\mathbf{y}_i = \mathbf{H}_i \mathbf{F}_i \mathbf{s}_i + \mathbf{H}_i \mathbf{F}_{\bar{i}} \mathbf{s}_{\bar{i}} + \mathbf{n}_i, \qquad \bar{i} = (N+1) - i. \quad (72)$$

Similarly, the channel $\mathbf{H}_i$ in (2) and the achievable rate in (7) remain valid where $\mathbf{\Omega}_{\bar{i}} \triangleq \sum_{\bar{i}=1, \bar{i} \neq i}^{N} \mathbf{H}_i \mathbf{F}_{\bar{i}} \mathbf{F}_{\bar{i}}^H \mathbf{H}_i^H + \sigma^2 \mathbf{I}_{N_{r_i}}$.

Hence, the max-min optimization can be expressed as

$$\max_{\boldsymbol{\theta}, \{\mathbf{F}_i\}_{i=1}^N} \min \left\{ \{R_i\}_{i=1}^N \right\} \tag{73a}$$

$$\text{s.t.} \quad \theta_n \in \mathcal{F} \quad \forall n \in \{1, \cdots, N_{ris}\} \tag{73b}$$

$$\sum_{i=1}^N \text{Tr}(\mathbf{F}_i \mathbf{F}_i^H) \leq P_{max} \tag{73c}$$

$$\Pr\left(\boldsymbol{\theta}, \{\mathbf{F}_i\}_{i=1}^N, \theta^{(ob)}\right) \leq \rho \quad \theta^{(ob)} \in \mathcal{A} \tag{73d}$$

where $\Pr\left(\boldsymbol{\theta}, \{\mathbf{F}_i\}_{i=1}^N, \theta^{(ob)}\right)$ is the power scattered by the RIS towards the direction $\theta^{(ob)}$, when the vector of reflection coefficients is set to $\boldsymbol{\theta}$, which is expressed as

$$\Pr(\boldsymbol{\theta}, \{\mathbf{F}_i\}_{i=1}^N, \theta^{(ob)}) = \sum_{i=1}^N \text{Tr}\left(\mathbf{G}^{(ob)} \text{diag}(\boldsymbol{\theta}) \mathbf{U} \mathbf{F}_i \mathbf{F}_i^H \mathbf{U}^H \right.$$
$$\left. \text{diag}(\boldsymbol{\theta})^H (\mathbf{G}^{(ob)})^H \right) \tag{74}$$

The same Arimoto-Blahut structure can be applied to reformulate the expression for the achievable rate as in (11). The reformulated expression of the achievable rate in (15) can be extended to

$$\mathbb{E}\left[\log_2\left(\frac{\mathcal{CN}(\mathbf{W}_i \mathbf{y}_i, \boldsymbol{\Sigma}_i)}{\mathcal{CN}(\mathbf{0}, \mathbf{I}_{N_{r_i}})}\right)\right]$$
$$= 2\,\text{Re}\left(\text{Tr}\left(\boldsymbol{\Sigma}_i^{-1} \mathbf{W}_i \mathbf{G}^{(i)} \text{diag}(\boldsymbol{\theta}) \mathbf{U} \mathbf{F}_i\right)\right) - \text{Tr}(\boldsymbol{\Sigma}_i^{-1})$$
$$- \text{Tr}\left(\mathbf{F}_i^H \mathbf{U}^H \text{diag}(\boldsymbol{\theta})^H \left(\mathbf{G}^{(i)}\right)^H \mathbf{W}_i^H \boldsymbol{\Sigma}_i^{-1} \right.$$
$$\left. \mathbf{W}_i \mathbf{G}^{(i)} \text{diag}(\boldsymbol{\theta}) \mathbf{U} \mathbf{F}_i\right)$$
$$- \sum_{\bar{i}=1, \bar{i} \neq i}^N \text{Tr}\left(\mathbf{F}_{\bar{i}}^H \mathbf{U}^H \text{diag}(\boldsymbol{\theta})^H \left(\mathbf{G}^{(i)}\right)^H \mathbf{W}_i^H \boldsymbol{\Sigma}_i^{-1}\right.$$
$$\left. \mathbf{W}_i \mathbf{G}^{(i)} \text{diag}(\boldsymbol{\theta}) \mathbf{U} \mathbf{F}_{\bar{i}}\right)$$
$$- \sigma^2 \text{Tr}\left(\mathbf{W}_i^H \boldsymbol{\Sigma}_i^{-1} \mathbf{W}_i\right) - N_{r_i} \log_2(\det(\boldsymbol{\Sigma}_i)) + N_{r_i} \tag{75}$$

where $\mathbf{W}_i$ and $\boldsymbol{\Sigma}_i$ are updated according to (12) and (13). We can also directly extend the expression of the mask constraint in (20) to

$$\Pr\left(\boldsymbol{\theta}, \{\mathbf{F}_i\}_{i=1}^N, \theta^{(ob)}\right)$$
$$= \sum_{i=1}^N \text{Tr}\left(\mathbf{G}^{(ob)} \text{diag}(\boldsymbol{\theta}) \mathbf{U} \mathbf{F}_i \mathbf{F}_i^H \mathbf{U}^H \text{diag}(\boldsymbol{\theta})^H (\mathbf{G}^{(ob)})^H\right)$$
$$= \sum_{i=1}^N \text{Tr}\left(\text{diag}(\boldsymbol{\theta})^H \left(\mathbf{G}^{(ob)}\right)^H \mathbf{G}^{(ob)} \text{diag}(\boldsymbol{\theta}) \mathbf{U} \mathbf{F}_i \mathbf{F}_i^H \mathbf{U}^H\right)$$
$$= \boldsymbol{\theta}^H \left(\mathbf{Q}^{(ob)} \odot \mathbf{T}_1 + \cdots + \mathbf{Q}^{(ob)} \odot \mathbf{T}_N\right) \boldsymbol{\theta} \tag{76}$$

### A. Update the RIS Phase Shift Vector $\boldsymbol{\theta}$

To update $\boldsymbol{\theta}$, we follow the same steps as in the dual-beam case. We optimize $\boldsymbol{\theta}$ keeping $\mathbf{F}_i$, $\mathbf{W}_i$ and $\boldsymbol{\Sigma}_i$, $i = 1, \ldots, N$, fixed. Following the same steps from (16)–(19), the sub-problem (22) can be extended to

$$\max_{\boldsymbol{\theta}} \min \left\{ -\boldsymbol{\theta}^H \mathbf{E}_1 \boldsymbol{\theta} + 2\,\text{Re}\left(\boldsymbol{\theta}^H \mathbf{b}_1\right) + c_1, \ldots, \right.$$
$$\left. -\boldsymbol{\theta}^H \mathbf{E}_N \boldsymbol{\theta} + 2\,\text{Re}\left(\boldsymbol{\theta}^H \mathbf{b}_N\right) + c_N \right\}$$
$$\text{s.t.} \quad \boldsymbol{\theta}^H \mathbf{I}_{N_{ris}}(:,n) \mathbf{I}_{N_{ris}}(:,n)^H \boldsymbol{\theta} = 1 \quad \forall n$$
$$\boldsymbol{\theta}^H (\mathbf{Q}^{(ob)} \odot \mathbf{T}_1 + \cdots + \mathbf{Q}^{(ob)} \odot \mathbf{T}_N) \boldsymbol{\theta} \leq \rho \tag{77}$$

where $\mathbf{E}_i = \mathbf{A}_i \odot \mathbf{B}_i + \mathbf{C}_i \odot \mathbf{D}_i$, $\mathbf{A}_i = \left(\mathbf{G}^{(i)}\right)^H \mathbf{W}_i^H \boldsymbol{\Sigma}_i^{-1} \mathbf{W}_i \mathbf{G}^{(i)}$, $\mathbf{B}_i = \left(\mathbf{U} \mathbf{F}_i \mathbf{F}_i^H \mathbf{U}^H\right)^T$, $\mathbf{C}_i = \left(\mathbf{G}^{(i)}\right)^H \mathbf{W}_i^H \boldsymbol{\Sigma}_i^{-1} \mathbf{W}_i \mathbf{G}^{(i)}$, $\mathbf{D}_i = \sum_{\bar{i}=1, \bar{i} \neq i}^N \left(\mathbf{U} \mathbf{F}_{\bar{i}} \mathbf{F}_{\bar{i}}^H \mathbf{U}^H\right)^T$ and $c_i = -\sigma^2 \text{Tr}(\mathbf{W}_i^H \boldsymbol{\Sigma}_i^{-1} \mathbf{W}_i) - N_{r_i} \log_2(\det(\boldsymbol{\Sigma}_i)) + N_{r_i} - \text{Tr}(\boldsymbol{\Sigma}_i^{-1})$.

Similar to the dual-beam case, the extended optimization problem (77) is a typical QCQP and can be solved using CVX.

### B. Update the Transmit Precoding $\mathbf{F_i}$

The optimization sub-problem (23) can be extended to

$$\max_{\{\mathbf{F}_i\}_{i=1}^N} \min\left\{ M_1(\mathbf{F}_1, \cdots, \mathbf{F}_N), \cdots, M_N(\mathbf{F}_1, \cdots, \mathbf{F}_N) \right\}$$
$$\text{s.t.} \quad \sum_{i=1}^N \text{Tr}\left(\mathbf{F}_i \mathbf{F}_i^H\right) \leq P_{max} \tag{78}$$

where $M_i(\mathbf{F}_1, \cdots, \mathbf{F}_N) = -\sum_{\bar{i}=1, \bar{i} \neq i}^N \text{Tr}\left(\mathbf{F}_{\bar{i}}^H \mathbf{J}_i \mathbf{F}_{\bar{i}}\right) + 2\,\text{Re}\left(\text{Tr}\left(\mathbf{F}_i^H \mathbf{K}_i^H\right)\right) + c_i - \text{Tr}\left(\mathbf{F}_i^H \mathbf{J}_i \mathbf{F}_i\right)$. The sub-problem in (78) is then decomposed by solving for each $\mathbf{F}_i$ separately. The optimization sub-problems (24) and (31) can be extended to the general multi-beam case as

$$\max_{\mathbf{F}_i} \min\left\{ -\text{Tr}\left(\mathbf{F}_i^H \mathbf{J}_i \mathbf{F}_i\right) + 2\,\text{Re}\left(\text{Tr}\left(\mathbf{F}_i^H \mathbf{K}_i^H\right)\right) - v_i, \right.$$
$$- \text{Tr}\left(\mathbf{F}_i^H \mathbf{J}_2 \mathbf{F}_i\right) + o_{i,2}, \ldots,$$
$$\left. - \text{Tr}\left(\mathbf{F}_i^H \mathbf{J}_N \mathbf{F}_i\right) + o_{i,N} \right\}$$
$$\text{s.t.} \quad \text{Tr}\left(\mathbf{F}_i \mathbf{F}_i^H\right) \leq P_{max}^{(i)} \tag{79}$$

where $v_i = \sum_{\bar{i}=1, \bar{i} \neq i}^N \text{Tr}\left(\mathbf{F}_{\bar{i}}^H \mathbf{J}_i \mathbf{F}_{\bar{i}}\right) - c_i$, $o_{i,\bar{i}} = -\sum_{\bar{i}=1, \bar{i} \neq i}^N \text{Tr}\left(\mathbf{F}_{\bar{i}}^H \mathbf{J}_{\bar{i}} \mathbf{F}_{\bar{i}}\right) + 2\,\text{Re}\left(\text{Tr}\left(\mathbf{F}_{\bar{i}}^H \mathbf{K}_{\bar{i}}^H\right)\right) + c_{\bar{i}}$. The optimization problem (79) is a typical QCQP, which is solved by CVX. Therefore, the problem can be generalized to the multi-beam beamforming case.

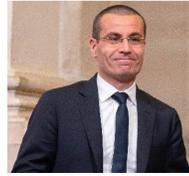

**Marco Di Renzo** (Fellow, IEEE) received the Laurea (cum laude) and Ph.D. degrees in electrical engineering from the University of L'Aquila, Italy, in 2003 and 2007, respectively, and the Habilitation à Diriger des Recherches (Doctor of Science) degree from University Paris-Sud (currently Paris-Saclay University), France, in 2013. Currently, he is a CNRS Research Director (Professor) and the Head of the Intelligent Physical Communications group with the Laboratory of Signals and Systems (L2S) at CNRS & CentraleSupélec, Paris-Saclay University, Paris, France, as well as a Chair Professor in Telecommunications Engineering with the Centre for Telecommunications Research – Department of Engineering, King's College London, London, United Kingdom. He was a France-Nokia Chair of Excellence in ICT at the University of Oulu (Finland), a Tan Chin Tuan Exchange Fellow in Engineering at Nanyang Technological University (Singapore), a Fulbright Fellow at The City University of New York (USA), a Nokia Foundation Visiting Professor at Aalto University (Finland), and a Royal Academy of Engineering Distinguished Visiting Fellow at Queen's University Belfast (U.K.). He is a Fellow of the IEEE, IET, EURASIP, and AAIA; an Academician of AIIA; an Ordinary Member of the European Academy of Sciences and Arts, an Ordinary Member of the Academia Europaea, and Ordinary Member of the Italian Academy of Technology and Engineering; an Ambassador of the European Association on Antennas and Propagation; and a Highly Cited Researcher. His recent research awards include the Michel Monpetit Prize conferred by the French Academy of Sciences, the IEEE Communications Society Heinrich Hertz Award, and the IEEE Communications Society Marconi Prize Paper Award in Wireless Communications. He served as the Editor-in-Chief of IEEE Communications Letters from 2019 to 2023. His current main roles within the IEEE Communications Society include serving as a Voting Member of the Fellow Evaluation Standing Committee, as the Chair of the Publications Misconduct Ad Hoc Committee, and as the Director of Journals. Also, he is on the Editorial Board of the Proceedings of the IEEE.

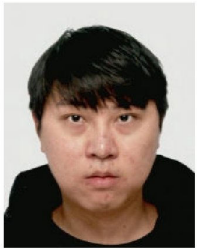

**Shumin Wang** is a postdoctoral researcher at Université Paris-Saclay, France, and a visiting researcher at Ranplan Wireless in Cambridge, UK. He received his Ph.D. in Networks, Information and Communication Sciences from Université Paris-Saclay in 2024. Prior to that, he earned his M.Sc. in Electronic Science and Technology from Central South University, Changsha, China, in 2021, and his B.Eng. in Engineering from Guangxi University, Nanning, China, in 2018. His research interests include Reconfigurable intelligent surfaces (RIS), optimization algorithms, and machine learning for wireless networks.

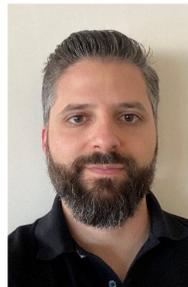

**Marios I. Poulakis** Dr.Eng. Marios I. Poulakis received the Diploma degree in Electrical and Computer Engineering (ECE) from the National Technical University of Athens (NTUA), Greece, and the M.Sc. degree in Management and Economics of Telecommunication Networks from the National and Kapodistrian University of Athens (NKUA), Greece, in July 2006 and December 2008, respectively. In May 2014, he received the Dr.Eng. degree from the NTUA. In 2007, he joined the Mobile Radio Communications Laboratory, ECE, NTUA as an associate researcher / project engineer participating in various industry and research-oriented projects. Moreover, he has teaching experience in undergraduate and postgraduate courses of ECE-NTUA and NKUA. From April 2016 to September 2017, M. Poulakis has been a post-doctoral researcher at the Department of Digital Systems, University of Piraeus. In addition, from April 2016 to October 2018, M. Poulakis has been an R&D Engineer at Feron Technologies P.C. Since November 2018, he has been a Wireless Technology Planning Engineer - Wireless and AI Expert at Huawei Technologies, Sweden. His research interests include wireless and satellite communications with emphasis on optimization mechanisms, machine learning and algorithm design. He has published more than 25 papers in international journals, conference proceedings and book chapters, while he has been a reviewer for international journals. He has several programming skills (such as C/C++, Java, Python and Android Dev). He is a Senior Member of IEEE.

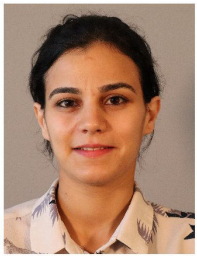

**Hajar El Hassani** is an associate professor at ENSEA and a researcher at ETIS Lab, CY Cergy Paris University, CNRS, France. She received her Ph.D. in electrical engineering from CY Cergy Paris University in 2022. Prior to her current role, she was a postdoctoral researcher at L2S, CNRS, CentraleSupélec, Université Paris-Saclay. She was also a visiting researcher at PUC-Rio in Rio de Janeiro, Brazil in 2022. She was awarded the Huawei "Seeds for the Future" scholarship in 2018. Her research focuses on modeling and optimization of wireless communication systems, with particular emphasis on green and energy-efficient technologies including Reconfigurable Intelligent Surfaces (RIS) and Ambient Backscatter Communications (AmBC), as well as resource allocation and machine learning for wireless networks. She is currently serving as the Grant Holder Scientific Representative of the COST Action 6G-PHYSEC.